\documentclass[12pt,a4paper,reqno]{amsart}
\usepackage[english]{babel}
\usepackage{amssymb,latexsym,amsfonts,amsthm,upref,amsmath}
\usepackage[foot]{amsaddr}
\usepackage[margin=1in]{geometry} 
\usepackage[dvipsnames,x11names]{xcolor}
 \definecolor{myblue}{HTML}{003399}
\usepackage{hyperref}
\hypersetup{colorlinks,citecolor=myblue,filecolor=black,linkcolor=myblue,urlcolor=myblue}
\usepackage{enumerate}  
\usepackage{tikz}
\usepackage{float}
\usepackage{cleveref}
\usepackage{mathtools}
\usepackage{cite}
\usepackage{filecontents}
\usepackage{enumitem}
\makeatletter
\newcommand{\leqnomode}{\tagsleft@true}
\newcommand{\reqnomode}{\tagsleft@false}

\makeatother
\newtheorem*{thm*}{Theorem}
\newtheorem*{lem*}{Lemma}
\newtheoremstyle{prim}{}{}{\normalfont}{}{\bfseries}{.}{ }{}
\newtheoremstyle{stil}{}{}{\slshape}{}{\bfseries}{.}{ }{}
\theoremstyle{stil}
\newtheorem{thm}{Theorem}[section]
\newtheoremstyle{defi}{}{}{}{}{\bfseries}{.}{ }{}
\theoremstyle{defi}
\newtheorem{defn}[thm]{Definition}
\theoremstyle{defi}
\newtheorem{rem}[thm]{Remark}
\theoremstyle{stil}
\newtheorem*{mthm*}{Main Theorem}
\newtheorem*{kor*}{Corollary}
\newtheorem{pro}[thm]{Proposition}
\theoremstyle{stil}
\newtheorem{lem}[thm]{Lemma}
\theoremstyle{stil}
\newtheorem{kor}[thm]{Corollary}
\theoremstyle{prim}

\newenvironment{prf}{\noindent \textit{Proof.}}{\null\hfill$\qed$\hskip
2mm\vskip 2mm}

\newcommand{\Yy}{ \text{\bfseries  Y}_{\hspace{-2pt}RTT}(\wvr{R})}
\newcommand{\Yp}{ \text{\bfseries  Y}_{\hspace{-2pt}RTT}^+(R)}
\newcommand{\Ypbar}{ \text{\bfseries  Y}_{\hspace{-2pt}RTT}^+(\wvr{R})}
\newcommand{\Ypc}{ \text{\bfseries  Y}_{\hspace{-2pt}RTT}^+(R)_c}
\newcommand{\Ybp}{ \text{\bfseries  Y}^+\hspace{-1pt}(R)}
\newcommand{\Ybpbar}{ \text{\bfseries  Y}^+\hspace{-1pt}(\wvr{R})}
\newcommand{\Ybpc}{ \text{\bfseries  Y}^+\hspace{-1pt}(R)_c}
\newcommand{\Yb}{ \text{\bfseries  Y}(R)}

\newcommand{\vac}{\mathop{\mathrm{\boldsymbol{1}}}}

\newcommand{\CC}{\mathbb{C}}

\newcommand{\Pc}{\mathcal{P}}

\newcommand{\Sc}{\mathcal{S}}

\newcommand{\Vcr}{ \mathcal{V}_c(R) }
\newcommand{\Mcr}{ \mathcal{M}(R)_c }
\newcommand{\Mcrphi}{ \mathcal{M}^{\phi}(R)_c }
\newcommand{\Mcrphir}{ \mathcal{M}^{\phi}_{RTT}(R)_c }

\newcommand{\R}{\wvr{R}}
\newcommand{\Rt}{\wht{R}}
\newcommand{\K}{\wvr{L}}
\newcommand{\T}{\wvr{S}}
\newcommand{\ttt}{\wvr{T}}
\newcommand{\tttt}{\wvr{t}}

\newcommand{\Dcrbar}{ \text{\bfseries  D}(\wvr{R})_c }

\newcommand{\gl}{\mathfrak{gl}}

\newcommand{\wht}{\widehat}
\newcommand{\wvr}{\overline}

\newcommand{\ot}{\otimes}
\newcommand{\ts}{\hspace{1pt}}
\newcommand{\qdet}{ {\rm qdet}\hspace{1pt}}
\newcommand{\tr}{ {\rm tr}}

\newcommand{\ndo}{\mathop{\mathrm{End}}}
\newcommand{\om}{\mathop{\mathrm{Hom}}}

\newcommand{\z}{\mathfrak z}

\newcommand{\cdotrl}{\mathop{\hspace{-2pt}\underset{\text{RL}}{\cdot}\hspace{-2pt}}}
\newcommand{\cdotlr}{\mathop{\hspace{-2pt}\underset{\text{LR}}{\cdot}\hspace{-2pt}}}

\newcommand{\fand}{\quad\text{and}\quad}
\newcommand{\Fand}{\qquad\text{and}\qquad}

\newcommand{\non}{\nonumber}
\newcommand{\beq}{\begin{equation}}
\newcommand{\eeq}{\end{equation}}
\newcommand{\ben}{\begin{equation*}}
\newcommand{\een}{\end{equation*}}

\makeatletter
\def\smalloverbrace#1{\mathop{\vbox{\m@th\ialign{##\crcr\noalign{\kern3\p@}%
  \tiny\downbracefill\crcr\noalign{\kern3\p@\nointerlineskip}%
  $\hfil\displaystyle{#1}\hfil$\crcr}}}\limits}
\makeatother

\makeatletter
\def\smallunderbrace#1{\mathop{\vtop{\m@th\ialign{##\crcr
   $\hfil\displaystyle{#1}\hfil$\crcr
   \noalign{\kern3\p@\nointerlineskip}%
   \tiny\upbracefill\crcr\noalign{\kern3\p@}}}}\limits}
\makeatother

\begin{document}

\title[ On the $h$-adic quantum vertex algebras associated   with Hecke  symmetries]
{On the $h$-adic quantum vertex algebras associated   with Hecke  symmetries}

\author{Slaven Ko\v{z}i\'{c}} 
\address{Department of Mathematics, Faculty of Science, University of Zagreb, Bijeni\v{c}ka cesta 30, 10\,000 Zagreb, Croatia}
\email{kslaven@math.hr}

\begin{abstract}
We study the quantum vertex algebraic framework for the Yangians of RTT-type and the braided Yangians associated with Hecke symmetries, 
introduced by Gurevich and Saponov. 
First, we construct several families of modules for the aforementioned Yangian-like algebras which, in the RTT-type case,   lead to a certain $h$-adic quantum vertex algebra $\mathcal{V}_c (R)$ via the Etingof--Kazhdan  construction, while, in the braided case, they produce ($\phi$-coordinated) $\mathcal{V}_c (R)$-modules. Next, we show that the coefficients of suitably defined quantum determinant  can be used to obtain central elements of   $\mathcal{V}_c (R)$, as well as the invariants of such ($\phi$-coordinated) $\mathcal{V}_c (R)$-modules. Finally, we investigate a certain algebra which is closely connected with the representation theory of $\mathcal{V}_c (R)$.
\end{abstract}

\maketitle

\allowdisplaybreaks


\section{Introduction}\label{intro}
\numberwithin{equation}{section}

The Yangians for classical Lie algebras present an important class of    quantum groups  which goes back to   Drinfeld \cite{Dri2}. 
Before the introduction of this family of  Hopf algebras, the algebra structure
of the Yangian  
 for the 
general linear Lie algebra 
$\mathfrak{gl}_N$ 
was  already intensively studied in the context of inverse scattering method; see, e.g., the papers by Takhtajan--Faddeev \cite{TF}, Kulish--Sklyanin \cite{KS}  and Tarasov \cite{T1,T2}.   
It is  a canonical deformation of the
universal enveloping algebra   of the  infinite-dimensional Lie algebra 
$\mathfrak{gl}_N [t]$.
Moreover, its    presentation can be given in terms of single ternary defining relation on the matrix generators, the so-called {\em RTT-relation},  which employs the Yang $R$-matrix.
 For more information on the Yangians  for classical Lie algebras  see the book by Molev \cite{Mol}.

In this paper, we consider  two classes of  generalizations of the Yangian for $\mathfrak{gl}_N$, which were introduced by Gurevich and Saponov \cite{GS}, the {\em Yangians of RTT-type} and the {\em braided Yangians}. 
 They are   defined in terms of   Yangian-like  RTT-relations, which, instead of the  Yang $R$-matrix, employ different $R$-matrices
obtained by the so-called Baxterization procedure
from the   involutive and Hecke symmetries. 
 Such algebras exhibit  properties similar to Yangians and have important   applications to Gaudin-type models; see \cite{GS,GSS}.
For simplicity, we consider only the (skew-invertible) Hecke symmetry case, although we briefly discuss the major differences  which occur for involutions as well.

The  goal of this paper is to associate the quantum vertex algebra  theory with the aforementioned   classes of Yangian-like algebras.
Our main motivation is the rich interplay, which appears in the classical theory, between representations of vertex algebras and infinite-dimensional Lie algebras, see, e.g., the books by E. Frenkel and Ben-Zvi \cite{FBZ}, I. Frenkel, Lepowsky and Meurman \cite{FLM} and Kac \cite{Kac2}.
Hopefully, in analogy with the classical case,   construction of such  new examples of quantum vertex algebras and  their applications
 might lead to  the better understanding of the general theory. For more information on the theory of quantum vertex algebras see, e.g., the papers by Etingof and Kazhdan \cite{EK5}, Li \cite{Li, Li1},  De Sole,   Gardini and Kac \cite{DGK} and references therein.

  Both the RTT-type and the braided Yangians can be associated with the $R$-matrix with additive and with  multiplicative spectral parameter (i.e., more precisely,  which satisfies the additive or the multiplicative version of the quantum Yang--Baxter equation). In all such cases,  we construct certain $\CC[[h]]$-modules generated by the quantum analogues of creation operators, such that they are equipped with a certain Yangian-like action   which resembles  the annihilation operators; see Prop. \ref{pro21},  \ref{pro32}, \ref{pro32bar} and \ref{proorp}.

In the case of Yangian of RTT-type associated with the additive $R$-matrix, this is essentially the Etingof--Kazhdan construction \cite{EK5}, so that we obtain a certain $h$-adic quantum vertex algebra $\Vcr$, where $c\in\CC$ and $R$ is the corresponding skew-invertible Hecke-symmetry; see Thm. \ref{mainthm1}.
However, in contrast with \cite{EK5}, the underlying braiding does no longer need to be of   type $1+\mathcal{O}(h)$, due to the different form of the $R$-matrix.
 On the other hand,   the case of braided Yangian for the additive $R$-matrix leads   to a family of $\Vcr$-modules; see Thm. \ref{glavni1}. 

As for the multiplicative $R$-matrix,  the actions of the corresponding RTT-type Yangian and   braided Yangian    produce   families of {\em $\phi$-coordinated modules} for $\Vcr$; see Thm. \ref{thm411} and \ref{mcrrr}. This is a certain   class of   structures in the   vertex algebra theory  introduced by Li \cite{Li1} in order to establish the connection with quantum groups; see also the more recent paper \cite{JKLT} for more information  on the theory of $\phi$-coordinated modules. Furthermore,  for the multiplicative RTT-type Yangian, we introduce in Subsect. \ref{subsec432} a certain new algebra $\Dcrbar$ via the defining relations resembling the RTT-presentation of the quantum affine algebra in type $A$. We show that a certain wide family of its modules is naturally equipped with the structure of $\phi$-coordinated 
$\Vcr$-module.  

In Sect. \ref{sec732}, motivated by the constructions from \cite{GS} and using the properties of the skew-symmetrizer established therein, we   study the {\em quantum determinant}, a certain formal power series in $\Vcr[[z]]$. We show that, under certain assumptions on the Hecke symmetry $R$, its coefficients belong to the center of the $h$-adic quantum vertex algebra $\Vcr$. Moreover, their images with respect to the ($\phi$-coordinated) module maps lead to the  ($\phi$-coordinated) module invariants. Finally, we apply such construction  to the algebra $\Dcrbar$, thus obtaining a family of  elements of its center.

\section{Preliminaries}\label{sec01}
\setcounter{equation}{0}

In this section,  we recall some basic properties of skew-invertible Hecke symmetries and their  $R$-matrices. All concepts are introduced over the commutative ring $\CC[[h]]$ so that they are compatible with the $h$-adic quantum vertex algebra theory. They are obtained from   the usual notions, defined  over the field $\CC(q)$, by setting $q=e^h$.
\subsection{Hecke symmetries} \label{subsec11}
 
Let $N\geqslant 2$ be an integer. Denote by $I$ and $P$ the identity and the permutation operator on $\CC^N \ot \CC^N$, respectively,
\beq\label{intr5}
I=\sum_{i,j=1}^N e_{ii}\ot e_{jj}\fand
P=\sum_{i,j=1}^N e_{ij}\ot e_{ji},
\eeq
where $e_{ij}$ are matrix units.   Let $h$ be a formal parameter and $\CC[[h]]$ the commutative ring of formal Taylor series  in $h$. Define a  {\em Hecke symmetry} over $\CC[[h]]$ as an element of $\ndo\CC^N \ot \ndo\CC^N [[h]]$ which satisfies the  {\em braid relation}
\beq\label{hecke1}
R_{12}\ts R_{23}\ts R_{12}
=
R_{23}\ts R_{12}\ts R_{23}
\eeq
and the condition
\beq\label{hecke2}
\left(R-e^h  I \right)\left(R+e^{-h}I \right)=0.
\eeq
In   the braid relation we use the standard tensor notation, where for any   $A=\sum A_1\ot A_2$ in $\ndo \CC^N \ot \ndo \CC^N [[h]]$ and distinct indices $1\leqslant r, s\leqslant m$ we write $A_{rs}$ for the element which acts as $A_1$ (resp. $A_2$) on the $r$-th (resp. $s$-th) tensor copy of $(\CC^N)^{\ot m}  $, i.e.   we have
$$
A_{rs}=\sum \left(1^{\ot (r-1)}\ot A_1 \ot 1^{\ot (m-r)}\right)
\left(1^{\ot (s-1)}\ot A_2 \ot 1^{\ot (m-s)}\right).
$$
In particular, in \eqref{hecke1} we have $m=3$ and $(r,s)\in\left\{(1,2),(2,3)\right\}$.
Regarding the identity in \eqref{hecke2},   the term $e^{\pm h}$ stands for the formal power series $\sum_{k\geqslant 0} (\pm h)^k/k!\in\CC[[h]]$.
Clearly, the identity implies that  $R$ is invertible and, furthermore, that its inverse is given by  $R^{-1}=R-(e^h -e^{-h})I$. Thus the constant term of $R$ with respect $h$, i.e. $R|_{h=0}$ is nonzero.

 The  Hecke symmetry  $R$ is said to be {\em skew-invertible} if there exists  an element $\Psi$ in $\ndo \CC^N \ot \ndo \CC^N [[h]]$ such that  
\beq\label{skew}
\tr_2\ts R_{12}\ts \Psi_{23}
=
P= \tr_2\ts  \Psi_{12}\ts R_{23},
\eeq
where the trace is taken over the second tensor factor.  As with $R$,  this implies that $\Psi|_{h=0}  $ is nonzero.
Throughout this paper, we consider only skew-invertible   Hecke symmetries. 
One well-known example of such   symmetry is given by  (cf. \cite{Dri,Jim}) 
$$
R =\sum_{i,j=1}^N e^{\delta_{ij} h} e_{ij}\ot e_{ji}+(e^h-e^{-h})\sum_{i<j}e_{ii}\ot e_{jj}.
$$
It can be easily checked that it satisfies the skew-invertibility condition \eqref{skew} with
$$
\Psi =\sum_{i,j=1}^N e^{-\delta_{ij}h}  e_{ij}\ot e_{ji}-(e^h-e^{-h})\sum_{i<j} e^{2(i-j)h}e_{ii}\ot e_{jj}.
$$ 
Observe that the evaluation   of $R$ at $h=0$   equals the permutation operator $P$, which is an example of  {\em skew-invertible involutive symmetry}, i.e. it is an involution which satisfies the braid relation   \eqref{hecke1} and the  skew-invertibility condition \eqref{skew} (with $\Psi = P$).

In this paper, we shall often use the ordered product notation, where the subscript of the product symbol determines the order of tensor factors. More precisely, for any   elements $A=\sum_i A_1^{(i)}\ot A_2^{(i)}$ and $B=\sum_j B_1^{(j)} \ot B_2^{(j)}$  in $\ndo \CC^N \ot \ndo \CC^N$ we write
$$
A\cdotlr B =\sum_{i,j} A_1^{(i)}B_1^{(j)}\ot  B_2^{(j)} A_2^{(i)}\Fand
A\cdotrl B =\sum_{i,j} B_1^{(j)}A_1^{(i)}\ot   A_2^{(i)}B_2^{(j)}.
$$
Such notation directly extends to multiple tensor factors, as well as to formal power series with coefficients in the tensor product algebra.
One easily checks  that the skew-invertibility condition \eqref{skew} is equivalent to   
\beq\label{rlprod}
(PR)\cdotrl (P\Psi) =I= (PR)\cdotlr (P\Psi)  .
\eeq
Conjugating the above equalities by the permutation operator $P$ we also find
\beq\label{rlprod2}
(RP)\cdotlr (\Psi P) =I= (RP)\cdotrl (\Psi P) .
\eeq
Finally, by applying the matrix transposition ${}^{t}\colon e_{ij}\mapsto e_{ji}$ to the $i$-th tensor factor of the  identities in \eqref{rlprod} and \eqref{rlprod2} with $i=1,2$ we get
\beq\label{trans}
(PR)^{t_i} (P\Psi)^{t_i} = (P\Psi)^{t_i} (PR)^{t_i}=I\fand
(RP)^{t_i} (\Psi P)^{t_i} = ( \Psi P)^{t_i} ( RP)^{t_i}=I .
\eeq
For  more information on     symmetries and their properties  we refer the reader to  \cite{GPS,O}.

\subsection{Normalizing series} \label{subsec1002}
Let us denote by  $\CC_* (x_1,\ldots ,x_n)$  
 the localization  of the ring of formal Taylor series
$\CC[[x_1,\ldots ,x_n]]$ at nonzero polynomials $\CC[x_1,\ldots ,x_n]^{\times}$.
There exists a unique embedding $ \CC_* (x_1,\ldots ,x_n)\to \CC((x_1))\ldots ((x_n))$. By extending  this embedding   to the $h$-adic completion of $\CC_* (x_1,\ldots ,x_n)$ we obtain the map
\beq\label{iotas}
\iota_{x_1,\ldots ,x_n}\colon \CC_* (x_1,\ldots ,x_n)[[h]]\to \CC((x_1))\ldots ((x_n))[[h]].
\eeq

Let $M>0$ be an integer.
 There exists  a unique   series $\wvr{g}(x)$ in $1+x\CC[[x,h]]$   satisfying 
\beq\label{snj6}
\wvr{g}(x)\ts \wvr{g}(xe^{-2h})\ldots \wvr{g}(xe^{-2(M-1)h})=\frac{1-xe^{-2(M-1)h}}{1-x};
\eeq
see \cite{FR} and \cite[Sect. 2]{KM} for more information.
Therefore, we have the identity
\beq\label{normal}
\wvr{f}(x)\ts \wvr{f}(xe^{-2h})\ldots \wvr{f}(xe^{-2(M-1)h})=\frac{1-xe^{-2(M-1)h}}{1-xe^{-2Mh}}\quad\text{for}\quad
\wvr{f}(x)=\frac{1-x}{1-xe^{-2h}}\ts\wvr{g}(x).
\eeq
The expressions in \eqref{snj6} and \eqref{normal} are understood as elements of
$\CC[[x,h]]\subset \CC((x))[[h]]$ via the embedding $\iota_x$, so that, in particular, we have
\beq\label{rcl6}
\iota_x\ts (1-xe^{-ah})^{-1} =\sum_{k\geqslant 0} x^ke^{-akh}\quad\text{for}\quad a\in\CC.
\eeq
As with the aforementioned expressions, throughout the paper  we usually omit the embedding symbol. In the multiple variable case,
we employ the   expansion convention where the choice of the embedding is determined by the order of the variables. More specifically, if $\sigma$ is a permutation in the symmetric group $\mathfrak{S}_n$, then
$(x_{\sigma_1}+\ldots + x_{\sigma_n})^r$ stands for   $\iota_{x_{\sigma_1},\ldots , x_{\sigma_n}}(x_{\sigma_1}+\ldots + x_{\sigma_n})^r$. For example, note that we have $(x_1+x_2)^r \neq (x_2+x_1)^r$ for $r<0$ as, by the convention, the former expression  should be expanded in the negative powers of $x_1$ and the latter in the negative powers of $x_2$.

Let $a\in\CC$ be nonzero.
As demonstrated in \cite[Sect. 2]{KM},    one can apply  the  substitution $x=e^{-2u/a}$ to $\wvr{f}(x)$ and then the embedding $\iota_u$, which produces the formal power series
$$
f(u)= \iota_u\ts  \wvr{f}(x)\big|_{x={e^{-2u/a}}}\in\CC((u))[[h]].
$$
The series $f(u)$ is invertible in $\CC((u))[[h]]$. Moreover, by \eqref{normal}  it satisfies
\beq\label{normal5}
f(u)\ts f(u+ah)\ldots f(u+a(M-1)h)= \frac{1-e^{-2u/a} e^{-2(M-1)h}}{1-e^{-2u/a}e^{-2M h}},
\eeq
where the right hand side is again understood as an element of $\CC((u))[[h]]$ via $\iota_u$.

The series $f(u)$ and $\wvr{f}(x)$ depend on the choice of the  integer $M$. In Section \ref{sec732}, we shall require that $M$ is equal to the rank of a certain Hecke symmetry. However,  the results of    other sections do not depend on the choice of $M$, so we suppress it in our notation.

\subsection{$R$-matrix with additive spectral parameter} \label{subsec12}
We  follow \cite[Sect. 4]{GS} to associate the $R$-matrix with a skew-invertible Hecke symmetry $R$. As before,  $a $ denotes an arbitrary nonzero complex number.
Consider the $R$-matrix $R(u)=R_{12}(u)$ defined by
\beq\label{rat2}
R(u)=\psi\ts f(u)\left(PR+\frac{e^h -e^{-h}}{e^{  2u/a} -1}P\right)\in\ndo\CC^N\ot\ndo\CC^N ((u))[[h]],
\eeq
where  $\psi$ is the invertible element of $\CC[[h]] $ such that $R(u)$ possesses the {\em unitarity property},
\beq\label{U}
R_{12}(u)\ts R_{21}(-u) =1.
\eeq
The constant $\psi$    can be found  as in the proof of \cite[Prop. 2.1]{KM}. 
The $R$-matrix \eqref{rat2}  satisfies the {\em Yang--Baxter equation},
\beq\label{YBE}
R_{12}(u)\ts R_{13}(u+v)\ts R_{23}(v)=
R_{23}(v)\ts R_{13}(u+v)\ts R_{12}(u).
\eeq
In addition to these properties, we shall often use the identities such as 
\begin{align}
R_{21}(u)\ts (RP)_{13}\ts (RP)_{23}
=(RP)_{23}\ts (RP)_{13}\ts R_{21}(u),\label{ybelike}\\
R_{31}(u)\ts (RP)_{23}\ts (RP)_{21}
=(RP)_{21}\ts (RP)_{23}\ts R_{31}(u).\label{ybelike2}
\end{align}
They can be  easily verified by using the braid relation \eqref{hecke1},  which is satisfied by $R$ and   $P$, along with  the explicit form  \eqref{rat2} of the $R$-matrix. Finally, note that $R(u)$ is obtained   from \cite[Prop. 12]{GS} by  multiplying the corresponding $R$-matrix by the permutation operator $P$ (from the left) and   the normalization term $\psi  f(u)$\footnote{Regarding the normalization, the term $f(u)$ in \eqref{rat2}, as well as $\wvr{f}(x)$ in \eqref{trig1} below, is used as it leads to nice properties of    quantum determinants, as we  demonstrate  in Subsections \ref{sec732b} and \ref{sec732c} below.}.

To simplify the following calculation, we  write
$$
f_1(u)=\psi\ts f(u)\fand f_2(u)=\frac{e^h -e^{-h}}{e^{  2u/a} -1},\qquad
\text{so that}\qquad R(u) = f_1(u) \left(PR + f_2 (u)P\right).
$$
The first family of identities in \eqref{trans} implies that the transposed $R$-matrices
$$
R(u)^{t_i}=f_1 (u) \left(PR\right)^{t_i}\left(I + f_2 (u)\left(P\Psi\right)^{t_i }P^{t_i}\right)\quad \text{with}\quad i=1,2 
$$
are invertible. Indeed, one easily checks that for
\beq\label{infsum}
S^{(i)}(u)= f_1(u)^{-1}\sum_{l\geqslant 0} \left(-f_2(u) \left(P\Psi\right)^{t_i }P^{t_i}\right)^{l}\left(P\Psi\right)^{t_i }
\in\ndo\CC^N\ot\ndo\CC^N ((u))[[h]] 
\eeq
we have
\beq\label{ccsym}
R(u)^{t_i}\ts S^{(i)}(u)=S^{(i)}(u)\ts R(u)^{t_i}=1\quad \text{for}\quad i=1,2.
\eeq
Notice that the infinite sum in \eqref{infsum} is convergent with respect to the $h$-adic topology since the series $f_2(u)$ belongs to $h  \CC  ((u))[[h]]$.   Let $S(u)=S^{(i)}(u)^{t_i}$.  The index $i$ is   redundant in this expression as $S^{(1)}(u)^{t_1}=S^{(2)}(u)^{t_2}$.
By applying the transposition on the $i$-th tensor factor of \eqref{ccsym} we rewrite the identities in terms of   ordered products as
\beq\label{ccsym2}
R(u) \cdotlr S (u) =R(u) \cdotrl S (u) =1 .
\eeq

\subsection{$R$-matrix with multiplicative spectral parameter}\label{subsec13} 

 Consider the $R$-matrix
\beq\label{trig1}
\R(x)=\wvr{f}(x)\left(PR+\frac{(e^h -e^{-h})x}{1-x}P\right)\in  \ndo\CC^N \ot    \ndo\CC^N [[x,h]];
\eeq
see \cite[Prop. 12]{GS}.
Roughly speaking, it is obtained from $\psi^{-1} R(u)$ by   replacing $e^{2u/a}$ by $x^{-1}$ and then applying the embedding $ \iota_x$; recall \eqref{iotas} and \eqref{rcl6}.  
More precisely, the connection between these $R$-matrices  is   established as follows. For any   $n> 0$ there exists $r\geqslant 0$ such that
$
(1-x)^r \R(x)
$
belongs to $\CC[x,h]$ modulo $h^n$. 
Indeed, the integer $r$ can be chosen so that  $(1-x)^r$ cancels all the poles of $\R(x)$ at $x=1$ modulo $h^n$.
Thus, we have  
\beq\label{conrmat}
\left( (1-x)^r \ts\R(x)\right)\big|_{x=e^{-2u/a}}^{\text{mod }h^n}
=
 \big(1-e^{-2u/a} \big)^r  \psi^{-1} R(u) \mod h^n .
\eeq
The superscript $\text{mod }h^n$ on the left hand side indicates that the expression inside the brackets is regarded modulo $h^n$, so that the substitution $x=e^{-2u/a}$ is indeed well-defined.

The $R$-matrix $\R(x)$ satisfies the multiplicative version  of the {\em Yang--Baxter equation},
\beq\label{YBEtrig}
\R_{12}(x)\ts \R_{13}(xy)\ts \R_{23}(y)=
\R_{23}(y)\ts \R_{13}(xy)\ts \R_{12}(x) 
\eeq
and   the identities \eqref{ybelike} and \eqref{ybelike2}.
Moreover,   arguing as in the end of Subsection \ref{subsec12}, one easily checks that there exists
 $\T (x)\in\ndo\CC^N \ot\ndo\CC^N [[x,h]]$ such that
\beq\label{ccsym2trig}
\R(x) \cdotlr \T (x) =\R(x) \cdotrl \T (x) =1 .
\eeq
As with   $\R(x)$, for any   $n> 0$ there exists $r\geqslant 0$  such that
\beq\label{conrmats}
\left( (1-x)^r \ts\T(x)\right)\big|_{x=e^{-2u/a}}^{\text{mod }h^n}
=
 \big(1-e^{-2u/a} \big)^r  \psi\ts S(u) \mod h^n .
\eeq

Observe that multiplying \eqref{rat2} and \eqref{trig1} by suitable factors we obtain  the $R$-matrices
$$
R'(u)=(1-e^{ - 2u/a})PR+(e^h -e^{-h})e^{ - 2u/a} P
\fand
\R'(x)=
(1-x)PR+ (e^h -e^{-h})x P
$$
which posses only nonnegative powers of $u$ and $x$, respectively. Furthermore, they satisfy
\beq\label{erpls}
 \R'(x) \big|_{x=e^{-2u/a}} 
=
R'(u) .
\eeq

Throughout the paper we shall  often write bar on the top of the  symbol to indicate that the corresponding object  comes with the multiplicative spectral parameter. In addition,  the letters   $u$, $v$ will usually indicate additive  and $x$, $y$ multiplicative variables.


\section{$h$-adic quantum vertex algebra  \texorpdfstring{$\Vcr$}{Vc(R)} }\label{sec02}
\setcounter{equation}{0}

In this section, 
we use the structure of the Yangian of RTT-type   to introduce   the creation     and   annihilation operators  over a certain $\CC[[h]]$-module $\Vcr$. By employing these operators we carry out the Etingof--Kazhdan-type construction    of the $h$-adic quantum vertex algebra   over $\Vcr$ associated with the skew-invertible Hecke symmetry $R$.  

The definition of the following algebra is motivated by the notion of Yangian of RTT-type, as given in \cite{GS}. We discuss a connection between these two classes of algebras in  Remark \ref{RTT-rem} below. Also, this algebra    can be regarded as the quantum function algebra associated to the $R$-matrix $R(u)$; cf.  \cite{EK4}.
Define
$\Yp$
as the topologically free associative algebra over the ring $\CC[[h]]$
generated by the elements $t_{ij}^{(-r)}$, where $i,j=1,\ldots ,N$ and $r=1,2,\ldots ,$ subject to the defining relations
\beq\label{rtt}
R(u-v)\ts T_1^+(u)\ts T_2^+(v)=T_2^+(v)\ts T_1^+(u) \ts R(u-v).
\eeq
The elements $T^+(u)$ are defined by
\beq\label{teplus}
T^+(u)=\sum_{i,j=1}^N e_{ij}\ot t_{ij}^+(u),\quad\text{where}\quad t_{ij}^+(u)=\sum_{r\geqslant 1} t_{ij}^{(-r)} u^{r-1},
\eeq
while the  subscripts in \eqref{rtt} indicate the   factors in the tensor product algebra, so that 
\beq\label{tftf}
T_1^+(u)=\sum_{i,j=1}^N e_{ij}\ot 1\ot  t_{ij}^+(u)\fand
T_2^+(v)=\sum_{i,j=1}^N 1\ot e_{ij}\ot   t_{ij}^+(v).
\eeq
In the above definition, we require that the given algebra is   topologically free as a $\CC[[h]]$-module; cf. \cite[Ch. XVI]{Kas}. This is easily accomplished by suitably modifying the ideal of its defining relations; see, e.g., \cite[Rem. 4.1]{c18} for more information. 

Let $V$ be a $\CC[[h]]$-module. We shall denote by $V((z))_h$ (resp. $V[z^{-1}]_h$) the $\CC[[h]]$-module of all formal power series $a(z)=\sum_{r\in\mathbb{Z}} a_r z^r$ in $V[[z^{\pm 1}]]$ (resp. in $V[[z^{- 1}]]$) such that we have $a_r \to 0$ when $r\to -\infty$ with respect to the $h$-adic topology. Such notation
naturally extends to the multiple variable case, so we write, for example, $V((z_1,\ldots ,z_n))_h$.
The following proposition can be proved by a direct calculation which relies on the properties  of the $R$-matrix  and the defining relations for the algebra $\Yp$; cf. \cite[Lemma 2.1]{EK5}.
\begin{pro}\label{pro21} 
For any $c\in \CC$ there exists a unique formal power series 
$ T^-(u)$ in
 $\ndo\CC^N \ot \om( \Yp,\Yp ((u))_h )$  
satisfying $T^-(u)\vac=I\ot \vac$  such that for all integers $n\geqslant 1$ we have
\begin{align}
&T^-_0(u)\ts T_1^+(v_1)\ldots T_n^+(v_n) = 
R_{01}(u-v_1+hc/2)^{-1}\ldots R_{0n}(u-v_n+hc/2)^{-1}\non\\
 &\qquad\qquad\qquad \times
T_1^+(v_1)\ldots T_n^+(v_n)\ts
R_{0n}(u-v_n-hc/2)\ldots R_{01}(u-v_1-hc/2).\label{Tu}
\end{align}
Moreover,   $T^-(u)$  is invertible in 
 $\ndo\CC^N \ot \om( \Yp,\Yp ((u))_h )$ and it 
   satisfies the following identities for operators on $\Yp$:
\begin{gather}
R(u-v )\ts T^-_1(u)\ts T^-_2 (v)=T^-_2 (v)\ts T^-_1(u)\ts R(u-v ),\label{rtt2} \\
R(u-v+hc/2)\ts T^-_1(u)\ts T_2^+(v)=T_2^+(v)\ts T^-_1(u)\ts R(u-v-hc/2) .\label{rtt3}
\end{gather}
\end{pro}

In order to emphasize that the $\CC[[h]]$-module $\Yp$ is regarded with respect to the  action of $T(u)$, which   depends on the choice of $c\in\CC$, we denote it by  $\Ypc$

\begin{rem}\label{RTT-rem}
The notion of Yangian of RTT-type    was introduced in \cite{GS} via defining relations of the form \eqref{rtt2}, 
$$
R(u-v )\ts T_1(u)\ts T_2 (v)=T_2 (v)\ts T_1(u)\ts R(u-v ),
$$
where the corresponding generator series $T(u)$ consists of nonpositive powers of the variable $u$. Regarding the operator series  $T^-(u)$, although satisfying the   relation of the same form, it possesses    nonnegative powers of $u$ as well, which appears to be common for the annihilation operators in Etingof--Kazhdan's construction; cf. \cite[Lemma 2.1]{EK5}.
\end{rem}

Let $u=(u_1,\ldots ,u_n)$ be a family of variables and $z$ a  single variable. We shall often use the following notation for   operators on $(\ndo\CC^N)^{\ot n} \ot\Yp$:
\beq\label{3not4}
T^{\pm}_{[n]}(u)=T^{\pm}_1(u_1)\ldots T^{\pm}_n(u_n)
\fand
T^{\pm}_{[n]}(z+u)=T^{\pm}_1(z+u_1)\ldots T^{\pm}_n(z+u_n),
\eeq
where $z+u=(z+u_1,\ldots ,z+u_n)$.
In addition, to simplify the notation, for any $b\in\CC$ we write $z+u+bh$ for the variables 
$ (z+u_1+bh,\ldots ,z+u_n+bh)$.
 Suppose $v=(v_1,\ldots ,v_m)$ is another family of variables. We introduce the following $R$-matrix products with entries in $(\ndo\CC^N)^{\ot n} \ot (\ndo\CC^N)^{\ot m}$:
\beq\label{r56}
R_{nm}(z+u-v+bh)=\prod_{1\leqslant i\leqslant n }^{\longrightarrow}\prod_{n+1\leqslant j\leqslant n+m}^{\longleftarrow} R_{ij}(z+u_i -v_{j-n}+bh),
\eeq
where the arrows indicate the order of factors. If the variable $z$ is omitted in the $R$-matrices on the right side, we denote the resulting expression by 
$R_{nm}(u-v+bh)$. Finally, we shall use the analogous notation for the products of the matrices $S(u)=S^{(i)}(u)^{t_i}$.
Using this notation, one can   generalize  \eqref{rtt}, \eqref{rtt2} and \eqref{rtt3} as follows; cf. \cite[Eq. (2.9)]{EK5}.
\begin{pro}\label{propositionz278}
For any $n,m\geqslant 1$ and the variables $u=(u_1,\ldots ,u_n)$, $v=(v_1,\ldots ,v_m)$   the following identities hold for operators on  $\Ypc$ :
\begin{gather*}
R_{nm}^{12}(u-v)\ts T_{[n]}^{\pm 13}(u)\ts   T_{[m]}^{\pm  23}(v)
=
T_{[m]}^{\pm  23}(v)\ts T_{[n]}^{\pm  13}(u)\ts R_{nm}^{12}(u-v),\\
R_{nm}^{12}(u-v+hc/2)\ts T_{[n]}^{- 13}(u)\ts   T_{[m]}^{+  23}(v)
=
T_{[m]}^{+ 23}(v)\ts T_{[n]}^{-  13}(u)\ts R_{nm}^{12}(u-v-hc/2).
\end{gather*}
\end{pro}

From now on, the tensor products of $\CC[[h]]$-modules are understood as
$h$-adically completed.
Essentially, the next theorem is  the Etingof--Kazhdan construction \cite[Thm. 2.3]{EK5} for the $R$-matrix \eqref{rat2}. As with the original version, it can be proved by directly     verifying the $h$-adic quantum vertex algebra axioms using the properties of the $R$-matrix  and the RTT-relations among the   operators $T^{\pm  }(u)$. We omit these calculations  since they go in parallel with  the case of   suitably normalized Yang $R$-matrix, whose details can be found in the proofs  of \cite[Thm. 2.3.8]{G} and \cite[Thm. 4.1]{JKMY}. Regarding the terminology,     by  the term {\em $h$-adic quantum vertex algebra}  we refer to  the generalization  of the notion of quantum VOA \cite[Sect. 1.4.1]{EK5}, as   given in \cite[Def. 2.20]{Li}.

\begin{thm}\label{mainthm1}
For any $c\in\CC$ there exists a unique structure of $h$-adic quantum vertex algebra on   $\Ypc$ such that the vacuum vector is the unit $\vac\in\Ypc$, the vertex operator map $Y(\cdot, z)$ is given by
\beq\label{vmap}
Y(T_{[n]}^+(u)\vac,z)= T_{[n]}^+(z+u)\ts T_{[n]}^-(z+u+hc/2)^{-1}
\eeq
and the braiding $\Sc(z)$ satisfies  the identity
\begin{align}
&\Sc^{34}(z)\left(R_{nm}^{12}(z+u-v)^{-1}\ts T_{[m]}^{+24}(v)\ts R_{nm}^{12}(z+u-v-hc)\ts T_{[n]}^{+13} (u)(\vac\ot\vac)\right)\non\\
&\qquad\qquad
=
T_{[n]}^{+13} (u)\ts R_{nm}^{12}(z+u-v+hc)^{-1}\ts  T_{[m]}^{+24}(v)\ts R_{nm}^{12}(z+u-v)(\vac\ot\vac).\label{braiding74}
\end{align}
\end{thm}

We  shall  denote the $h$-adic quantum vertex algebra established in Theorem \ref{mainthm1} by $\Vcr$. Regarding the formula for the braiding,  note that the superscripts $1,2,3,4$ in \eqref{braiding74} indicate the tensor factors as follows:
\beq\label{subskripti}
\smalloverbrace{(\ndo\CC^N)^{\ot n}}^{1} \ot \smalloverbrace{(\ndo\CC^N)^{\ot m}}^{2}\ot \smalloverbrace{\Vcr}^{3} \ot  \smalloverbrace{\Vcr}^{4} .
\eeq

\begin{rem}\label{skiisrem}
The results of this section,  including Theorem \ref{mainthm1}, also hold  if $R$ is a  skew-invertible involutive symmetry. In this case, one can employ, e.g.,    the $R$-matrix    
\beq\label{skiis}
R(u)= \left(1-\frac{ah}{u}\right)^{-1}\left(PR +\frac{ah}{u}P\right)\in\ndo\CC^N\ot\ndo\CC^N [[h/u]] 
\eeq
from \cite[Sect. 4]{GS} instead of \eqref{rat2}.
Thus the corresponding  operator series $ T^-(u)$, defined as in Proposition \ref{pro21}, possesses only nonpositive powers of $u$, so that it  produces the action of the  Yangian of  RTT-type  associated with $R$ \cite{GS} (defined over the ring $\CC[[h]]$); cf. Remark \ref{RTT-rem}.
In particular, if $R=P$, by suitably rescaling the generators   one obtains from    $T^-(u)$    the  action of the Yangian for $\gl_N$     on the corresponding  dual Yangian, which comes from the structure of the level $c$ double Yangian for $\gl_N$; see \cite{I,Kho}.
\end{rem}

\section{\texorpdfstring{$\Vcr$}{Vc(R)}-modules}\label{subsec31}

In this section, we use the structure of the braided Yangian associated to the $R$-matrix \eqref{rat2}   to construct a family of $\Vcr$-modules.

\subsection{Braided Yangian-type action (additive case)}

The definition of the following algebra is motivated by the notion of braided Yangian    from \cite{GS}. A connection between these two classes of algebras is discussed  in  Remark \ref{RTTplus-rem} below. 
Let
$\Ybp$
be the topologically free associative algebra over the ring $\CC[[h]]$
generated by the elements $l_{ij}^{(-r)}$, where $i,j=1,\ldots ,N$ and $r=1,2,\ldots ,$ subject to the defining relations
\beq\label{rttbraided}
R_{12}(u-v)\ts L_1^+(u)\ts (RP)_{12}\ts L_2^+(v)=L_2^+(v)\ts (RP)_{21}\ts L_1^+(u) \ts R_{21}(u-v).
\eeq
The elements $L^+(u)$ are defined by
\beq\label{lplus}
L^+(u)=\sum_{i,j=1}^N e_{ij}\ot l_{ij}^+(u),\quad\text{where}\quad l_{ij}^+(u)=\delta_{ij}-\sum_{r\geqslant 1} l_{ij}^{(-r)} u^{r-1},
\eeq
while their subscripts indicate the   factors in the tensor product algebra as in \eqref{tftf}.  

Let $u=(u_1,\ldots ,u_n)$ be a family of variables and $z$ a single variable. We introduce the following notation for the formal power series with coefficients in $(\ndo\CC^N)^{\ot n}\ot \Ybp$:
\beq\label{561}
L_{[n]}^+ (z+u)=\prod_{1\leqslant i\leqslant n }^{\longrightarrow}
L_i^+(z+u_i)\ts (RP)_{i\ts i+1}\ldots (RP)_{i\ts n}
.
\eeq
If  $z$ is omitted from the right hand side, we denote the resulting expression    by $L_{[n]}^+ (u)$.

\begin{lem}\label{lemma31}
The coefficients of all matrix entries of $L_{[n]}^+(u)=L_{[n]}^+(u_1,\ldots ,u_n)$ with $n\geqslant 1$ along with the unit $1$ span an $h$-adically dense $\CC[[h]]$-submodule of $\Ybp$.
\end{lem}

\begin{prf}
The   lemma follows from the identity
\beq\label{562}
L_{[n]}^+ (u)=\prod_{1\leqslant i\leqslant n }^{\longrightarrow}
L_i^+(u_i)\ts (RP)_{i\ts i+1}\ldots (RP)_{i\ts n}.
\eeq
Indeed, it is clear that the coefficients of all matrix entries of $L_{1}^+(u_1)\ldots L_n^+(u_n)$ with $n\geqslant 1$ along with   $\vac$ span an $h$-adically dense $\CC[[h]]$-submodule of $\Ybp$. However,   due to \eqref{rlprod2}, we can   remove all terms $RP$ from the right hand side of \eqref{562}   by multiplying the equality by suitable copies of $\Psi P$. Thus, we can express $L_{1}^+(u_1)\ldots L_n^+(u_n)$   in terms of $L_{[n]}^+(u_1,\ldots ,u_n)$, which implies the lemma.
\end{prf}

We have the following braided Yangian analogue of Proposition \ref{pro21}:

\begin{pro}\label{pro32} 
For any $c\in \CC$ there exists a unique formal power series 
$ L^-(u)$ in
 $\ndo\CC^N \ot \om( \Ybp,\Ybp ((u))_h )$  
satisfying $L^-(u)1=I\ot 1$  such that for all integers $n\geqslant 1$ and the variables  $v=(v_1,\ldots , v_n)$  we have
\begin{align}
 &L^{-0}(u)\ts (RP)_{01}\ldots (RP)_{0n}\ts L_{[n]}^{+12}(v)  = 
R_{01}(u-v_1+hc/2)^{-1}\ldots R_{0n}(u-v_n+hc/2)^{-1}\non
\\
&\qquad\qquad\times
L_{[n]}^{+12}(v) \ts (RP)_{10}\ldots (RP)_{n0}\ts
R_{n0}(u-v_n-hc/2)\ldots R_{10}(u-v_1-hc/2).\label{Lu}
\end{align}
 Moreover, $L^-(u)$ is invertible  in $\ndo\CC^N \ot \om( \Ybp,\Ybp ((u))_h )$   
and it   satisfies the following identities for operators on $\Ybp$:
\begin{align}
R_{12}(u-v)\ts L_1^-(u)\ts (RP)_{12}\ts L_2^-(v) &=L_2^-(v)\ts (RP)_{21}\ts L_1^-(u) \ts R_{21}(u-v),\label{ll2} \\
R_{12}(u-v+hc/2)\ts L_1^-(u)\ts (RP)_{12}\ts L_2^+(v)&=L_2^+(v)\ts (RP)_{21}\ts L_1^-(u) \ts R_{21}(u-v-hc/2) .\label{rll3}
\end{align}
\end{pro}

 \begin{prf}
First, note that  the superscripts in \eqref{Lu} indicate tensor factors  as follows:
$$
\smalloverbrace{(\ndo\CC^N) }^{0} \ot \smalloverbrace{(\ndo\CC^N)^{\ot n}}^{1}\ot \smalloverbrace{\Ybp}^{2}.
$$
The fact that $L^-(u)$ is well-defined by \eqref{Lu} can be proved by a direct calculation which shows that it preserves the ideal of defining relations \eqref{rttbraided} for the algebra $\Ybp$. The calculation relies on the Yang--Baxter equation \eqref{YBE} and the properties \eqref{ybelike}, \eqref{ybelike2} of the $R$-matrix. Next, by Lemma \ref{lemma31}  the expression \eqref{Lu} uniquely determines   $L^-(u)$. Regarding the invertibility, one easily checks that the inverse of $L^-(u)$ is given by
\begin{align*}
     L^{-0}(u)^{-1} 
L_{[n]}^{+12}(v)    
 = \big(   S_{n0} (-u+v_n-hc/2) \ldots S_{10} (-u+v_1-hc/2)\big) 
\cdotrl \Big(\hspace{-1pt} (RP)_{01}\ldots (RP)_{0n}  \Big.& \\ 
 \quad\Big. 
 \times L_{[n]}^{+12}(v)\ts R_{10}(u-v_1-hc/2)^{-1}\ldots R_{n0}(u-v_n-hc/2)^{-1}
(RP)_{n0}^{-1}\ldots (RP)_{10}^{-1} \Big).& 
\end{align*}
Hence both $L^- (u)$ and its inverse belong to $\ndo\CC^N \ot \om( \Ybp,\Ybp ((u))_h )$  due to the form of $R(u)$ and $S(u)$. Finally, the identities \eqref{ll2} and \eqref{rll3} are   immediate  consequences of \eqref{Lu}, which can be verified using   the properties \eqref{YBE}, \eqref{ybelike} and \eqref{ybelike2} of the $R$-matrix.
\end{prf}

\begin{rem}\label{RTTplus-rem}
The notion of braided Yangian  was introduced in \cite{GS} via defining relations  
$$
P_{12}\ts R_{12}(u-v )\ts L_1(u)\ts R_{12}\ts  L_1 (v)=L_1 (v)\ts R_{12}\ts L_1(u)\ts P_{12}\ts R_{12}(u-v ),
$$
where the corresponding generator series $L(u)$ consists of nonpositive powers of the variable $u$ and its constant term equals the identity. This relation is equivalent to
$$
R_{12}(u-v)\ts L_1 (u)\ts (RP)_{12}\ts L_2 (v)  =L_2 (v)\ts (RP)_{21}\ts L_1 (u) \ts R_{21}(u-v),
$$
i.e. it takes the form of \eqref{ll2}. However, the operators  $ L^-(u)$ considered in \eqref{ll2} contain both positive and negative powers of $u$.
\end{rem}

The action of $L^-(u)$ given by \eqref{Lu} depends on the choice of $c\in\CC$.
Hence, to indicate  that the $\CC[[h]]$-module $\Ybp$ is regarded with respect to such
action, we denote it from now on by $\Ybpc$. Generalizing the notation, we write $L_{[n]}^-(z+u)$ (resp. $L_{[n]}^-(u)$) for  the expression obtained from \eqref{561} (resp. \eqref{562}) by replacing all $L_i^+(z+u_i)$ (resp. $L_i^+(u_i)$) on the right hand side by the operator series $L_i^-(z+u_i)$ (resp. $L_i^-(u_i)$). Also, the original series in \eqref{561} and \eqref{562} are from now on regarded as operators on $\Ybpc$ with respect to the algebra multiplication.
 
The relations   established in Proposition \ref{pro32} can be generalized as follows. First,   denote by $R_{\bar{n}\bar{m}}(u-v+bh)$
the expression obtained from $R_{nm}(u-v+bh)$, where $b\in\CC$, by conjugating all $R$-matrices by the permutation operator $P$, i.e.
\beq\label{asin7}
R_{\bar{n}\bar{m}}( u-v+bh)=\prod_{1\leqslant i\leqslant n }^{\longrightarrow}\prod_{n+1\leqslant j\leqslant n+m}^{\longleftarrow} R_{ji}( u_i -v_{j-n}+bh).
\eeq
Also, we write $R_{\bar{n}\bar{m}}( z+u-v+bh)$ for the analogous product of the $R$-matrices $R_{ji}( z+u_i -v_{j-n}+bh)$; recall
\eqref{r56}. 
Furthermore, we use the analogous notation for the products of matrices $S(u)$. Thus, for example, by combining  \eqref{U} and
\eqref{ccsym2} one easily proves  
\beq\label{hjf1}
 R_{nm}(u-v)^{-1}\cdotlr S_{\bar{n}\bar{m}}(-u+v)=   R_{nm}(u-v)^{-1}\cdotrl S_{\bar{n}\bar{m}}(-u+v) =1.
\eeq
Next, we introduce the following notation
$$
(RP)_{nm}=\prod_{1\leqslant i\leqslant n }^{\longrightarrow}\prod_{n+1\leqslant j\leqslant n+m}^{\longrightarrow} (RP)_{ij}\Fand
(RP)_{\bar{n}\bar{m}}=\prod_{1\leqslant i\leqslant n }^{\longrightarrow}\prod_{n+1\leqslant j\leqslant n+m}^{\longrightarrow} (RP)_{ji}.
$$
The following proposition can be verified by a direct calculation which relies on the  relations \eqref{rttbraided}, \eqref{ll2} and \eqref{rll3} and the $R$-matrix identities \eqref{ybelike} and \eqref{ybelike2}.

\begin{pro}\label{propositionz}
For any $n,m\geqslant 1$ and the variables $u=(u_1,\ldots ,u_n)$, $v=(v_1,\ldots ,v_m)$   the following identities hold for operators on $\Ybpc$:
\begin{gather*}
R_{nm}^{12}(u-v)\ts L_{[n]}^{\pm 13}(u)\ts (RP)_{nm}^{12}\ts L_{[m]}^{\pm 23}(v)
=
L_{[m]}^{\pm 23}(v)\ts (RP)_{\bar{n}\bar{m}}^{12}\ts L_{[n]}^{\pm 13}(u)\ts R_{\bar{n}\bar{m}}^{12}(u-v),\\
R_{nm}^{12}(u-v+hc/2)\ts L_{[n]}^{- 13}(u)\ts (RP)_{nm}^{12}\ts L_{[m]}^{+ 23}(v)
=
L_{[m]}^{+ 23}(v)\ts (RP)_{\bar{n}\bar{m}}^{12}\ts L_{[n]}^{- 13}(u)\ts R_{\bar{n}\bar{m}}^{12}(u-v-hc/2).
\end{gather*}
\end{pro}

Note that by using \eqref{hjf1} we can rewrite the second identity in Proposition \ref{propositionz}   as
\begin{align}
&\,L_{[n]}^{- 13}(u)^{-1}   
L_{[m]}^{+ 23}(v)    =  S_{\bar{n}\bar{m}}^{12}(-u+v-hc/2)\non \\
 &\qquad\qquad\cdotrl\left(
 (RP)_{nm}^{12}\ts L_{[m]}^{+ 23}(v) \ts R_{\bar{n}\bar{m}}^{12}(u-v-hc/2)^{-1}\ts L_{[n]}^{- 13}(u)^{-1}\ts  ((RP)_{\bar{n}\bar{m}}^{12})^{-1}\right).\label{linvertz}
\end{align}
 
\begin{lem}\label{learestricted}
For any $c\in\CC$, $n\geqslant 1$ and the variables $u=(u_1,\ldots ,u_n)$ we have
$$
L_{[n]}^-(u)^{- 1} \in (\ndo\CC^N)^{\ot n}\ot \om( \Ybp_c,\Ybp_c ((u_1,\ldots,u_n))_h ).
$$ 
\end{lem}

\noindent\textit{Proof.}
By Lemma \ref{lemma31}, the coefficients of   matrix entries of all $L_{[m]}^+(v)1$ with $m\geqslant 1$ and $v=(v_1,\ldots ,v_m)$ along with $1$ span an $h$-adically dense $\CC[[h]]$-submodule of $\Ybp_c$. Thus it is sufficient to check that their images under 
$L_{[n]}^-(u)^{- 1}$ belong to 
$(\ndo\CC^N)^{\ot n}\ot  \Ybp_c ((u_1,\ldots,u_n))_h$.  This follows immediately    by applying   \eqref{linvertz} to $1\in \Ybp_c$ since $L_{[n]}^-(u)^{- 1} 1$ is constant with respect to the variables $u_1,\ldots ,u_n$. Indeed, due to the form of $R(u)$ and $S(u)$, as given by \eqref{rat2} and \eqref{infsum}, and the expansion convention from Subsection \ref{subsec1002}, the remaining terms on the right hand side,   $S_{\bar{n}\bar{m}}^{12}(-u+v-hc/2)$ and $R_{\bar{n}\bar{m}}^{12}(u-v-hc/2)^{-1}$ belong to
$$
(\ndo\CC^N)^{\ot n}\ot (\ndo\CC^N)^{\ot m} ((u_1,\ldots,u_n))_h [[v_1,\ldots ,v_m]]. \eqno\qed
$$ \vskip 2mm

\subsection{The    \texorpdfstring{$\Vcr$}{Vc(R)}-module structure over \texorpdfstring{$\Ybp_c$}{Y+(R)c}}\label{subsecx32}

In the next theorem we construct an example of $\Vcr$-module over the $\CC[[h]]$-module of $\Ybp_c$. 
For readers convenience, we first recall the definition of module for   $h$-adic quantum vertex algebra  \cite[Def. 2.23]{Li}.

\begin{defn}\label{mod}
Let $(\mathcal{V},Y,\vac, \Sc)$ be an $h$-adic quantum vertex algebra.
A {\em $\mathcal{V}$-module} is a pair $(\mathcal{M},Y_\mathcal{M})$ such that $\mathcal{M}$ is a topologically free $\CC[[h]]$-module and  
\begin{align*}
Y_\mathcal{M}(z)=Y_\mathcal{M}(\cdot, z) \colon  \mathcal{V} \ot \mathcal{M}&\to \mathcal{M}((z))_h\\
u\ot w&\mapsto Y_\mathcal{M}(z)(u\ot w)=Y_\mathcal{M}(u,z)w=\sum_{r\in\mathbb{Z}} u_r w \ts z^{-r-1} 
\end{align*}
is  a $\mathbb{C}[[h]]$-module map
which satisfies
 $Y_\mathcal{M}(\vac,z)w=w$  for all $w\in \mathcal{M}$ and  
  the {\em weak associativity}: for any $u,v\in  \mathcal{V} $, $w\in\mathcal{M}$ and integer $k\in\mathbb{Z}_{\geqslant 0}$ there exists $r\in\mathbb{Z}_{\geqslant 0}$ such that
\beq\label{wassoc7}
(z_0 +z_2)^r\ts Y_{\mathcal{M}}(u,z_0+z_2)\ts Y_{\mathcal{M}}(v, z_2)w
-(z_0 +z_2)^r\ts Y_{\mathcal{M}}(Y(u,z_0)v,z_2)w
\eeq
belongs to $h^k \mathcal{M} [[z_0^{\pm 1},z_2^{\pm 1}]]$, i.e. such that the above expressions coincide modulo $h^k$.
\end{defn}

To emphasize that $\Ybp_c$ is now regarded as a $\Vcr$-module, we denote it by $\Mcr$.   

\begin{thm}\label{glavni1}
For any $c\in\CC$ there exists a unique structure of module for the $h$-adic quantum vertex algebra $\Vcr$ on $\Mcr= \Ybp_c$ such that 
$Y_{\Mcr}(\vac,z)=1_{\Mcr}$ and such that
for all integers $n\geqslant 1 $ and the variables $u=(u_1,\ldots ,u_n)$ we have
\beq\label{modmap}
Y_{\Mcr}(T_{[n]}^+(u)\vac,z)
=
L_{[n]}^+(z+u)\ts
L_{[n]}^-(z+u+hc/2)^{-1}.
\eeq
\end{thm}

\begin{prf}
Let us show that the $\CC[[h]]$-module map is  well-defined by \eqref{modmap}, i.e. that it annihilates the ideal of defining relations \eqref{rtt}. Given the definition of the algebra $\Yp $, it is sufficient to check that for any indices $1\leqslant i <n $ and the  families of variables 
$$u=(u_1,\ldots ,u_n)\fand u^{(i)}=(u_1,\ldots, u_{i-1},u_{i+1},u_i,u_{i+2},\ldots ,u_n)$$ the image of  the expression
$$
R_{i\ts i+1}(u_i-u_{i+1}) \ts T_{[n]}^+(u)\vac
- P_{i\ts i+1}\ts T_{[n]}^+(u^{(i)})\ts P_{i\ts i+1} \vac \ts R_{i\ts i+1}(u_i-u_{i+1})   
$$
under $Y_{\Mcr}(\cdot,z)$  is trivial. In other words, due to \eqref{modmap}, we have to prove that
$$
R'_{i\ts i+1}   L_{[n]}^+(z+u) 
L_{[n]}^-(z+u+hc/2)^{-1}
-
 P_{i\ts i+1}  L_{[n]}^+(z+u^{(i)}) 
L_{[n]}^-(z+u^{(i)}+hc/2)^{-1}  P_{i\ts i+1} 
R'_{i\ts i+1}
$$
is zero for $R'_{i\ts i+1}=R_{i\ts i+1}(u_i-u_{i+1})$. However, this immediately follows from the identities 
$$
R'_{i\ts i+1}   L_{[n]}^{\pm}(z+u+bh) = P_{i\ts i+1}\ts  L_{[n]}^{\pm}(z+u^{(i)}+bh) \ts P_{i\ts i+1}\ts R_{i+1\ts i} (u_i-u_{i+1}) \quad\text{with}\quad b\in\CC,
$$
which can be verified by using    \eqref{ybelike2}  
and the
relations \eqref{rttbraided} and \eqref{ll2}.

Next, as the coefficients of   matrix entries of all $T_{[n]}^+(u)\vac$ with $n\geqslant 1$ along with $\vac$ span an $h$-adically dense $\CC[[h]]$-submodule of $\Vcr$,  the $\CC[[h]]$-module map $Y_{\Mcr}(\cdot,z)$ is uniquely determined by  \eqref{modmap}. Moreover, by Lemma \ref{learestricted}, its image   belongs to $\Mcr((z))_h$. 
Therefore, it remains to prove that it satisfies the   weak associativity, as given in    Definition \ref{mod}.

We start the proof of the weak associativity property by considering the image of
\beq\label{org4}
T_{[n]}^{+13}(u)\ts T_{[m]}^{+24}(v) (\vac\ot\vac)\in(\ndo\CC^N)^{\ot (n+m)} \ot\Vcr^{\ot 2}[[u_1,\ldots,u_n,v_1,\ldots ,v_m]] 
\eeq
under the   operators from \eqref{wassoc7}. 
Throughout the calculation, we use superscripts $1,2,3,4$ to indicate the 
tensor factors as in \eqref{subskripti}. 
We shall need the following simple consequence of the
  second relation in Proposition \ref{propositionz278}:
\begin{align}
T_{[n]}^{-  13}(u)^{-1}\ts T_{[m]}^{+ 23}(v) \vac
= &\,\,
  S_{\bar{n}\bar{m}}^{12}(-u+v-hc/2) \cdotrl\left(
   T_{[m]}^{+ 23}(v)\vac  \ts R_{nm}^{12}(u-v-hc/2)^{-1}\hspace{-1pt} \right).\label{linvertz2}
\end{align}

Due to \eqref{modmap}, the application of $Y_{\Mcr}(\cdot , z_0+z_2)(1\ot Y_{\Mcr}(\cdot ,z_2))$, which corresponds to 
the first term in \eqref{wassoc7},  to the expression in \eqref{org4} produces
\beq\label{cmpr1}
L_{[n]}^{+13}(z_0+z_2+u)\ts
L_{[n]}^{-13}(z_0+z_2+u+hc/2)^{-1}
L_{[m]}^{+23}( z_2+v)\ts
L_{[m]}^{-23}( z_2+v+hc/2)^{-1}.
\eeq
On the other hand, by applying $Y_{\Mcr}(\cdot, z_2)(Y(\cdot,z_0)\ot 1)$, which corresponds to 
the second term in \eqref{wassoc7},
  to \eqref{org4} and then using \eqref{vmap} and \eqref{linvertz2}
we get
\begin{align}
&S_{\bar{n}\bar{m}}^{12}(-z_0-u+v-hc)\cdotrl\left(
L_{[n+m]}^{+ 123 }(z_2+z_0+u,z_2+v)\right.\non\\
&\left.\times   L_{[n+m]}^{-  123}(z_2+z_0+u+hc/2,z_2+v+hc/2)^{-1}\ts R_{nm}^{12}(z_0+u-v)^{-1} \right),\label{cmpr2}
\end{align}
where $(z_2+z_0+u+bh,z_2+v+bh)$ with $b\in \CC$ denotes the $(n+m)$-tuple of variables
$$
(z_2+z_0+u_1+bh,\ldots ,z_2+z_0+u_n+bh,z_2+v_1+bh,\ldots , z_2+v_m+bh).
$$

We now compare the expressions in \eqref{cmpr1} and \eqref{cmpr2}. First, we use \eqref{linvertz} to swap the order of the two middle factors in \eqref{cmpr1}, thus getting
\begin{align*}
&
S_{\bar{n}\bar{m}}^{12}(-z_0-u+v-hc)\cdotrl\Big(
L_{[n]}^{+13}(z_0+z_2+u) 
(RP)_{nm}^{12}\ts L_{[m]}^{+ 23}(z_2+v)   \\
&\times R_{\bar{n}\bar{m}}^{12}(z_0+u-v)^{-1}\ts L_{[n]}^{-13}(z_0+z_2+u+hc/2)^{-1} \ts ((RP)_{\bar{n}\bar{m}}^{12})^{-1} L_{[m]}^{-23}( z_2+v+hc/2)^{-1}\Big).
\end{align*}
Next,   using the first family of relations in Proposition \ref{propositionz}, we reorder the last four factors:
\begin{align*}
&
S_{\bar{n}\bar{m}}^{12}(-z_0-u+v-hc)\cdotrl\Big(
L_{[n]}^{+13}(z_0+z_2+u) 
(RP)_{nm}^{12}\ts L_{[m]}^{+ 23}(z_2+v)   \\
&\times 
L_{[m]}^{-23}( z_2+v+hc/2)^{-1}\ts
((RP)_{nm}^{12})^{-1} \ts
L_{[n]}^{-13}(z_0+z_2+u+hc/2)^{-1}\ts
R_{nm}^{12}(z_0+u-v)^{-1}    \Big).
\end{align*}
Finally, it is clear that this equals to
\begin{align}
&S_{\bar{n}\bar{m}}^{12}(-z_0-u+v-hc)\cdotrl\left(
L_{[n+m]}^{+ 123 }(z_0+z_2+u,z_2+v)\right.\non\\
&\left.\times   L_{[n+m]}^{-  123}(z_0+z_2+u+hc/2,z_2+v+hc/2)^{-1}\ts R_{nm}^{12}(z_0+u-v)^{-1} \right).\label{cmpr3}
\end{align}

Note that the expressions in \eqref{cmpr2} and \eqref{cmpr3} are not equal as the term $L_{[n+m]}^{-  123}(\ldots)^{-1}$ in the former should be expanded in negative powers of $z_2$ and in the latter in negative powers of $z_0$.
Let $w\in\Mcr$ and let $a_1,\ldots, a_n,b_1,\ldots ,b_m,k $ be arbitrary positive integers. We shall now consider the coefficients of the monomials
\beq\label{monomi}
u_1^{a'_1}\ldots u_n^{a'_n} v_1^{b'_1}\ldots v_m^{b'_m}h^{k'}\quad\text{for}\quad
0\leqslant a'_i\leqslant a_i,\quad
0\leqslant b'_j\leqslant b_j,\quad
0\leqslant k' < k.
\eeq
By Lemma \ref{learestricted}, there exists   $r\geqslant 1$ such that the coefficients of monomials \eqref{monomi} in
\begin{align*}
&(z_0+z_2)^r\ts 
L_{[n+m]}^{-  123}(z_0+z_2+u+hc/2,z_2+v+hc/2)^{-1}w\Fand\\
&(z_0+z_2)^r\ts 
L_{[n+m]}^{-  123}(z_2+z_0+u+hc/2,z_2+v+hc/2)^{-1}w 
\end{align*}
coincide.
Therefore, suppose we apply    \eqref{cmpr2} and \eqref{cmpr3} to the vector $w$ and then  multiply  them by $(z_0+z_2)^r$. Then the  coefficients of the monomials \eqref{monomi} in the resulting expressions  coincide,
so we conclude that the weak associativity holds.
\end{prf}

\begin{rem}
As with Section \ref{sec02}, the results of this section can be also recovered when $R$ is a skew-invertible involutive symmetry  by using the $R$-matrix \eqref{skiis}; recall Remark \ref{skiisrem}.
In this case, the corresponding   series $ L^-(u)$ from Proposition \ref{pro32} 
possesses only nonpositive powers of $u$ and its constant term is the identity. Thus  it   
produces the action of the  braided Yangian associated with $R$   \cite{GS}; recall Remark \ref{RTTplus-rem}.  Naturally, we assume here that the braided Yangian is regarded over the ring  $\CC[[h]]$. 
\end{rem}

\section{ \texorpdfstring{$\phi$}{phi}-coordinated \texorpdfstring{$\Vcr$}{Vc(R)}-modules}\label{subsec32}

  In this section, we introduce certain modules for the braided Yangian and the RTT-type Yangian     associated with the $R$-matrix \eqref{trig1} defined over $\CC[[h]]$. Using the modules we construct two families of $\phi$-coordinated $\Vcr$-modules for $\phi(z_2,z_0)=z_2 e^{-2z_0 /a} $. In the end, we introduce a certain algebra   and establish a partial connection between its restricted modules  and  $\phi$-coordinated  $\Vcr$-modules.

\subsection{Action of the braided Yangian (multiplicative case)}\label{sbsc41}
The following  construction goes in parallel with Section \ref{subsec31}. 
We start with the topologically free associative algebra $\Ybpbar$ over the ring $\CC[[h]]$
which is 
generated by the elements $\wvr{l}_{ij}^{(-r)}$, where $i,j=1,\ldots ,N$ and $r=1,2,\ldots .$ These are subject to the defining relations
\beq\label{rttbraidedtrig}
\R_{12}(x/y)\ts \K_1^+(x)\ts (RP)_{12}\ts \K_2^+(y)=\K_2^+(y)\ts (RP)_{21}\ts \K_1^+(x) \ts \R_{21}(x/y).
\eeq
The elements $\K^+(x)$ are defined by
$$
\K^+(x)=\sum_{i,j=1}^N e_{ij}\ot \wvr{l}_{ij}^+(x),\quad\text{where}\quad \wvr{l}_{ij}^+(x)=\delta_{ij}-\sum_{r\geqslant 1} \wvr{l}_{ij}^{(-r)} x^{r-1}.
$$

Let $x=(x_1,\ldots ,x_n)$ be a family of variables. We introduce the following notation for the formal power series with coefficients in $(\ndo\CC^N)^{\ot n}\ot \Ybpbar$:
\beq\label{673}
\K_{[n]}^+ (x)=\prod_{1\leqslant i\leqslant n }^{\longrightarrow}
\K_i^+(x_i)\ts (RP)_{i\ts i+1}\ldots (RP)_{i\ts n}
.
\eeq
The next lemma can be proved analogously to   Lemma \ref{lemma31}.
\begin{lem}\label{sbt6}
The coefficients of all matrix entries of $\ts\K_{[n]}^+(x)=\K_{[n]}^+(x_1,\ldots ,x_n)$ with $n\geqslant 1$ along with the unit $1$ span an $h$-adically dense $\CC[[h]]$-submodule of $\Ybpbar$.
\end{lem}

We omit the proof of the following proposition as it goes in parallel with the proof of Proposition \ref{pro32}. In particular, it relies on Lemma \ref{sbt6} and  $R$-matrix properties from Subsection \ref{subsec13}.  The superscripts in the formula \eqref{Lubar} below are again used to indicate tensor factors analogously to the statement of Proposition \ref{pro32}. 

\begin{pro}\label{pro32bar} 
For any $c\in \CC$ there exists a unique formal power series 
$ \K^-(x)$ in
 $\ndo\CC^N \ot \om( \Ybpbar,\Ybpbar [x^{-1}]_h )$  
satisfying $\K^-(x)\vac=I\ot \vac$  such that for all integers $n\geqslant 1$  and a family of variables $y=(y_1,\ldots , y_n)$ we have
\begin{align}
 &\K^{-0}(x)\ts (RP)_{01}\ldots (RP)_{0n}\ts \K_{[n]}^{+12}(y)  = 
\R_{10}(y_1 e^{hc/a}/x ) \ldots \R_{n0}(y_n e^{hc/a}/x) \non
\\
&\qquad\qquad\times
\K_{[n]}^{+12}(y) \ts (RP)_{10}\ldots (RP)_{n0}\ts
\R_{0n}(y_n e^{ -hc/a} /x)^{-1}\ldots \R_{01}(y_1 e^{- hc/a} /x )^{-1}.\label{Lubar}
\end{align}
 The series $ \K^-(x)$ is invertible 
in  $\ndo\CC^N \ot \om( \Ybpbar,\Ybpbar [x^{-1}]_h )$  and it satisfies the following identities for   operators on $\Ybpbar$:
\begin{align}
\R_{12}(x/y)\ts \K_1^-(x)\ts (RP)_{12}\ts \K_2^-(y) &=\K_2^-(y)\ts (RP)_{21}\ts \K_1^-(x) \ts \R_{21}(x/y),\label{ll2bar} \\
 \K_1^-(x)\ts (RP)_{12}\ts \K_2^+(y)\ts \R_{12}( ye^{-hc/a}/x)&=\R_{21}(ye^{hc/a} /x )\ts\K_2^+(y)\ts (RP)_{21}\ts \K_1^-(x)  .\label{rll3bar}
\end{align}
\end{pro}

As before, we write $\Ybpbar_c$ to indicate that the $\CC[[h]]$-module
$\Ybpbar$ is regarded  with respect to the action \eqref{Lubar} of the series $\K(x)$, which depends on the choice of $c\in\CC$.
Also, generalizing the notation \eqref{673}, we denote by $\K_{[n]}^-(x)$   the expression obtained from \eqref{673} by replacing all   $\K_i^+(x_i)$ on the right hand side  by   $L_i^-(x_i)$. The original series in \eqref{673} are from now on regarded as operators on $\Ybpbar_c$ with respect to the algebra multiplication.

\begin{rem}
Note that \eqref{ll2bar} is equivalent to the defining relations of the braided Yangian $\Yb $  \cite{GS} associated with the $R$-matrix $\R(x)$. As the constant term of $\R(x) $ with respect to   $x$ equals $ PR$, one easily checks that the constant term of the   series $\K^-(x)$, as given by \eqref{Lubar}, is the identity. 
Suppose that $\Yb $ is defined over $\CC[[h]]$ and that its matrix of generators is denoted by $\K(x)$; cf. \cite[Sect. 6]{GS}.
Then for any $c\in\CC$ the assignment  $\K(x)\mapsto \K^-(x) $ defines the structure of $\Yb $-module over $\Ybpbar_c$.
\end{rem}

The commutation relations for operators $\K^{\pm 1}(x)$ can be generalized as follows. First,  for the families of variables $x=(x_1,\ldots,x_n)$, $y=(y_1,\ldots ,y_m)$ and $b\in\CC$ we write
$$
\R_{nm}( xe^{bh} /y)=\prod_{1\leqslant i\leqslant n }^{\longrightarrow}\prod_{n+1\leqslant j\leqslant n+m}^{\longleftarrow} \R_{ij}( x_i e^{bh}/y_{j-n}).
$$
As in \eqref{asin7}, we shall write bar on the top of the subscripts $n$ and $m$ to indicate that the $R$-matrices are conjugated by the permutation operator $P$,  so that we have
\beq\label{bsr5}
\R_{\bar{n}\bar{m}}( xe^{bh} /y)
=\prod_{1\leqslant i\leqslant n }^{\longrightarrow}\prod_{n+1\leqslant j\leqslant n+m}^{\longleftarrow} \R_{ji}( x_i e^{bh}/y_{j-n}).
\eeq
We   write "opp" in the superscript if the $R$-matrix factors appear in the opposite order, 
$$
\R_{nm}^{\text{opp}}( xe^{bh} /y)=\prod_{1\leqslant i\leqslant n }^{\longleftarrow}\prod_{n+1\leqslant j\leqslant n+m}^{\longrightarrow} \R_{ij}( x_i e^{bh}/y_{j-n}).
$$
In addition, we apply analogously the above notation to the products of matrices $\T(x)$; recall \eqref{ccsym2trig}. 
The next proposition can be verified by a direct calculation which relies on the 
   relations \eqref{rttbraidedtrig}, \eqref{ll2bar} and \eqref{rll3bar} and the $R$-matrix properties   from Subsection \ref{subsec13}. 

\begin{pro}\label{propositionz2}
For any $n,m\geqslant 1$ and the variables $x=(x_1,\ldots ,x_n)$, $y=(y_1,\ldots ,y_m)$   the following identities hold for operators on $\Ybpbar_c$:
\begin{gather*}
\R_{nm}^{12}(x/y)\ts \K_{[n]}^{\pm 13}(x)\ts (RP)_{nm}^{12}\ts \K_{[m]}^{\pm 23}(y)
=
\K_{[m]}^{\pm 23}(y)\ts (RP)_{\bar{n}\bar{m}}^{12}\ts \K_{[n]}^{\pm 13}(x)\ts \R_{\bar{n}\bar{m}}^{12}(x/y),\\
 \K_{[n]}^{- 13}(x)\ts (RP)_{nm}^{12}\ts \K_{[m]}^{+ 23}(y)\ts \R^{{\text{\normalfont opp}}\ts 12}_{nm}(y e^{ -hc/a}/x)
=
\R^{{\text{\normalfont opp}}\ts 12}_{\bar{n}\bar{m}}(y e^{ hc/a}/x)\ts
\K_{[m]}^{+ 23}(y)\ts (RP)_{\bar{n}\bar{m}}^{12}\ts \K_{[n]}^{- 13}(x) .
\end{gather*}
\end{pro}

\subsection{Action of the RTT-type Yangian (multiplicative case)}\label{sbsc42}
We proceed in parallel with Subsection \ref{sbsc41} to construct the action of the RTT-type Yangian. 
Let $\Ypbar$ be the topologically free associative algebra  over the ring $\CC[[h]]$
generated by the elements $\tttt_{ij}^{(-r)}$, where $i,j=1,\ldots ,N$ and $r=1,2,\ldots .$ These are subject to the defining relations
\beq\label{rrt0}
\R (x/y)\ts \ttt_1^+(x)\ts   \ttt_2^+(y)=\ttt_2^+(y) \ts \ttt_1^+(x) \ts \R (x/y).
\eeq
The elements $\ttt^+(x)$ are defined by
\beq\label{ghj76}
\ttt^+(x)=\sum_{i,j=1}^N e_{ij}\ot \tttt_{ij}^+(x),\quad\text{where}\quad \tttt_{ij}^+(x)=  \sum_{r\geqslant 1} \tttt_{ij}^{(-r)} x^{r-1}.
\eeq

 Consider   the formal power series
$\ttt_{[n]}^+ (x)\coloneqq  
\ttt_1^+(x_1) \ldots \ttt_n^+(x_n)
$
    whose  coefficients belong to $(\ndo\CC^N)^{\ot n}\ot \Ypbar$.
Clearly, the coefficients of all matrix entries of $\ts\ttt_{[n]}^+(x) $ with $n\geqslant 1$ along with the unit $1$ span an $h$-adically dense $\CC[[h]]$-submodule of $\Ypbar$.

We now give the RTT-type   counterpart of Proposition \ref{pro32bar}. Roughly speaking, it can be regarded as the multiplicative version of Proposition \ref{pro21}.
 It can  be easily proved using the properties  of the $R$-matrix  and the defining relations for the algebra  $\Ypbar$.

\begin{pro}\label{proorp} 
For any $c\in \CC$ there exists a unique formal power series 
$ \ttt^-(x)$ in
 $\ndo\CC^N \ot \om( \Ypbar,\Ypbar [x^{-1}]_h )$  
satisfying $\ttt^-(x)\vac=I\ot \vac$  such that for all integers $n\geqslant 1$ and the variables $y=(y_1,\ldots , y_n)$ we have
\begin{align}
 \ttt^{-0}(x)\ts   \ttt_{[n]}^{+12}(y)  = \,\,&
\R_{10}(y_1 e^{hc/a}/x ) \ldots \R_{n0}(y_n e^{hc/a}/x) \non
\\
& \times
\ttt_{[n]}^{+12}(y)  \ts
\R_{n0}(y_n e^{ -hc/a} /x)^{-1}\ldots \R_{10}(y_1 e^{- hc/a} /x )^{-1}.\label{tubar}
\end{align}
   The series $ \ttt^-(x)$ is invertible in  $\ndo\CC^N \ot \om( \Ypbar,\Ypbar [x^{-1}]_h )$  and
it satisfies the following identities for operators on $\Ypbar$:
\begin{align}
\R_{12}(x/y)\ts \ttt_1^-(x) \ts \ttt_2^-(y) &=\ttt_2^-(y)\ts  \ttt_1^-(x) \ts \R_{12}(x/y), \label{rrt} \\
 \ttt_1^-(x)\ts  \ttt_2^+(y)\ts \R_{21}( ye^{-hc/a}/x)&=\R_{21}(ye^{hc/a} /x )\ts\ttt_2^+(y)\ts   \ttt_1^-(x)  .\label{rrt2} 
\end{align}
\end{pro}

The commutation relations from the   proposition  can  be   easily generalized to the multiple operator  case, as was done for the braided Yangian setting in 
Proposition \ref{propositionz2}.
 
Again, we    use the subscript $c$ and write $\Ypbar_c$ to indicate that the $\CC[[h]]$-module
$\Ypbar$ is   regarded  with respect to the action \eqref{tubar} of   $\ttt(x)$, which depends on the choice of  $c\in\CC$.
Also,  we write $
\ttt_{[n]}^- (x)= 
\ttt_1^-(x_1) \ldots \ttt_n^-(x_n)
$. The   series $\ttt_{[n]}^{+ } (x)$  are from now on regarded as operators on $\Ypbar_c$ with respect to the algebra multiplication.

\begin{rem}
The relations \eqref{rrt} are equivalent to the defining relations of the RTT-type Yangian $\Yy$  \cite{GS} associated    \eqref{trig1}.
Suppose that $\Yy $ is defined over  $\CC[[h]]$ and that its matrix of generators is denoted by $\wvr{T}(x)$; cf. \cite[Sect. 6]{GS}.
Then for any $c\in\CC$ the assignment $\wvr{T}(x)\mapsto \ttt^-(x) $ defines
the structure of $\Yy $-module over $\Ypbar_c$.
\end{rem}

\subsection{Constructing  \texorpdfstring{$\phi$}{phi}-coordinated \texorpdfstring{$\Vcr$}{Vc(R)}-modules}\label{subsecx42}

The notion of $\phi$-coordinated module, where $\phi=\phi(z_2,z_0)$ is an  associate of the one-dimensional additive formal group,  was introduced by Li \cite{Li1}. In this paper  we consider its straightforward generalization over the ring $\CC[[h]]$; cf.\cite[Def. 5.1]{c18}. Regarding the associate, we set
\beq\label{assc5}
\phi(z_2,z_0)=z_2 e^{-2z_0 /a} 
.\eeq 
To present the precise definition, we need some notation. 
Let $V$ be a  topologically free $\CC[[h]]$-module.
Suppose we have  some elements
$A,C\in\om(V,V[[z_1^{\pm 1},z_2^{\pm 1}]])$ and $B\in\om(V,V((z_1,z_2))_h)$  satisfying  $A=B+h^k C$. To indicate that $A$ possesses such decomposition, we shall   write 
$A\in\om(V,V((z_1,z_2)))\mod h^k$. Note that the substitution
$B\left|_{z_1 =\phi(z_2,z_0)}\right.$
is well-defined although the substitution $A\left|_{z_1 =\phi(z_2,z_0)}\right.$
does not need to exist. We shall often denote $B\left|_{z_1 =\phi(z_2,z_0)}\right.$ by $A\big|_{z_1 =\phi(z_2,z_0)}^{\text{mod }h^k}\big.$; see, e.g.,   \eqref{associativitymod0} and \eqref{associativitymod} below.

\begin{defn}\label{phimod}
Let $(\mathcal{V},Y,\vac, \Sc)$ be an $h$-adic quantum vertex algebra.
A {\em $\phi$-coordinated $\mathcal{V}$-module} is a pair $(\mathcal{M}^\phi,Y_{\mathcal{M}^\phi})$ such that $\mathcal{M}^\phi$ is a topologically free $\CC[[h]]$-module and  
\begin{align*}
Y_{\mathcal{M}^\phi}(z)=Y_{\mathcal{M}^\phi}(\cdot, z) \colon  \mathcal{V} \ot \mathcal{M}^\phi&\to \mathcal{M}^\phi((z))_h\\
u\ot w&\mapsto Y_{\mathcal{M}^\phi}(z)(u\ot w)=Y_{\mathcal{M}^\phi}(u,z)w=\sum_{r\in\mathbb{Z}} u_r w \ts z^{-r-1} 
\end{align*}
is  a $\mathbb{C}[[h]]$-module map
which satisfies
 $Y_{\mathcal{M}^\phi}(\vac,z)w=w$  for all $w\in \mathcal{M}^\phi$ and  
  the {\em weak associativity}: for any $u,v\in  \mathcal{V} $, $k\in\mathbb{Z}_{\geqslant 0}$ there exists $r\in\mathbb{Z}_{\geqslant 0}$ such that
\begin{align}
&(z_1-z_2)^r\ts Y_{\mathcal{M}^\phi}(u,z_1)Y_{\mathcal{M}^\phi}(v,z_2)\in\om (\mathcal{M}^\phi,\mathcal{M}^\phi((z_1,z_2)) )\mod h^k  ,\label{associativitymod0}\\
&\big((z_1-z_2)^r\ts Y_{\mathcal{M}^\phi} (u,z_1)Y_{\mathcal{M}^\phi}(v,z_2)\big)\big|_{z_1= \phi(z_2,z_0)}^{\text{mod } h^k}  \big. \non\\
&\qquad-  (\phi(z_2,z_0) - z_2)^r\ts Y_{\mathcal{M}^\phi} (Y(u,z_0)v,z_2 )\ts
\in\ts  h^k \om(\mathcal{M}^\phi,\mathcal{M}^\phi[[z_0^{\pm 1},z_2^{\pm 1}]]).\label{associativitymod}
\end{align}
\end{defn}

In the next theorem, we construct an example of $\phi$-coordinated $\Vcr$-module using the structure of  braided Yangian module from Subsection \ref{sbsc41}. 
To emphasize that $\Ybpbar
_c$ is now regarded as a $\phi$-coordinated $\Vcr$-module, we denote it by $\Mcrphi$. 

\begin{thm}\label{thm411}
For any $c\in\CC$ there exists a unique structure of $\phi$-coordinated $\Vcr$-module    on $\Mcrphi= \Ybpbar_c$ such that for all integers $n\geqslant 1 $ and the families of variables $u=(u_1,\ldots ,u_n)$ and $x=(x_1,\ldots ,x_n)$ we have
\beq\label{modmap2}
Y_{\Mcrphi}(T_{[n]}^+(u)\vac,z)
=
\K_{[n]}^+(x)\big|_{x_i=\phi(z,u_i)}\ts
\K_{[n]}^-(x)^{-1}\big|_{x_i=\phi(z,u_i+hc/2)}.
\eeq
\end{thm}

\begin{prf}
We omit   some   technical  details, as they go in parallel with the proof of Theorem \ref{glavni1},
  and only discuss the differences which occur in this case.
First, regarding the expression in \eqref{modmap2}, 
 the substitutions $x_i=\phi(z,u_i)$ and $x_i=\phi(z,u_i+hc/2)$ are  assumed to be carried out simultaneously for all $i=1,\ldots ,n$.
Note that the
second family of relations in Proposition \ref{propositionz2} implies the identity
\begin{align}
&\K_{[n]}^{- 13}(x)^{-1}\ts 
\K_{[m]}^{+ 23}(y)  
=  
\T^{  12}_{\bar{n}\bar{m}}(y e^{ hc/a}/x)\non\\
&\qquad\qquad \cdotrl\left(
  (RP)_{nm}^{12}\ts \K_{[m]}^{+ 23}(y)  \ts \R^{{\text{\normalfont opp}}\ts 12}_{nm}(y e^{ -hc/a}/x)\ts \K_{[n]}^{- 13}(x)^{-1}\ts ((RP)_{\bar{n}\bar{m}}^{12})^{-1} 
	\right)
 .\label{prff54}
\end{align}
By applying \eqref{prff54} to $1 $ and then  arguing as in the proof of Lemma  \ref{learestricted} we find  
$$
\K_{[n]}^-(x)^{- 1} \in (\ndo\CC^N)^{\ot n}\ot \om( \Mcrphi,\Mcrphi [x_1^{-1},\ldots,x_n^{-1}]_h ),
$$
so that  the substitutions in \eqref{modmap2} are well-defined.

 The fact that the $\CC[[h]]$-module map $Y_{\Mcrphi}(\cdot,z)$ is well-defined by \eqref{modmap2} is verified as in the proof of Theorem \ref{glavni1}. However, the corresponding argument also employs   \eqref{erpls} and the following observation: if we replace the $R$-matrix $R(u-v)$ (resp. $\R(x/y)$) in     \eqref{rtt} (resp. \eqref{rttbraidedtrig} and \eqref{ll2bar}) by  $R'(u-v)$ (resp. $\R'(x/y)$), we obtain an equivalent equality.

Regarding \eqref{associativitymod0}, note that this requirement is   necessary in order for the substitution $z_1= \phi(z_2,z_0)$ in the first term of \eqref{associativitymod} to be well-defined. It is verified in a similar way to the equality between
\eqref{cmpr1} and \eqref{cmpr3}. Nonetheless, we give some details in order to take care of the differences. First, by applying 
$Y_{\mathcal{M}^\phi}(\cdot,z_1)(1\ot Y_{\mathcal{M}^\phi}(\cdot,z_2))$ to
\eqref{org4}
we obtain
$$
\K_{[n]}^{+13}(x)\big|_{x_i=\phi(z_1,u_i)}\ts
\K_{[n]}^{-13}(x)^{-1}\big|_{y_i=\phi(z_1,u_i+hc/2)}
\K_{[m]}^{+23}(y)\big|_{y_j=\phi(z_2,v_j)}\ts
\K_{[m]}^{-23}(y)^{-1}\big|_{y_j=\phi(z_2,v_j+hc/2)},
$$
where $x=(x_1,\ldots ,x_n)$ and $y=(y_1,\ldots ,y_m)$. Next, we   employ \eqref{prff54} to  swap the second and the third  term above, thus getting 
\begin{align}
\T^{  12}_{\bar{n}\bar{m}} 
\cdotrl & \left(
\K_{[n]}^{+13}(x)\big|_{x_i=\phi(z_1,u_i)}\ts
  (RP)_{nm}^{12}\ts \K_{[m]}^{+ 23}(y)\big|_{y_j=\phi(z_2,v_j)}   
	\right.
\non\\
  & \left. \times \R^{{\text{\normalfont opp}}\ts 12}_{nm}\ts  \K_{[n]}^{- 13}(x)^{-1}\big|_{x_i=\phi(z_1,u_i+hc/2)}\ts ((RP)_{\bar{n}\bar{m}}^{12})^{-1} \ts
	\K_{[m]}^{-23}(y)^{-1}\big|_{y_j=\phi(z_2,v_j+hc/2) }
	\right)
,\label{last4}
\end{align}
where
\begin{align*}
\T^{  12}_{\bar{n}\bar{m}}=
\T^{  12}_{\bar{n}\bar{m}}(y e^{ hc/a}/x)\Big|_{\substack{x_i=\phi(z_1,u_i+hc/2)\\ \hspace{-23pt}y_j=\phi(z_2,v_j) } }\fand \R^{{\text{\normalfont opp}}\ts 12}_{nm}=\R^{{\text{\normalfont opp}}\ts 12}_{nm}(ye^{-hc/a}  /x)\Big|_{\substack{x_i=\phi(z_1,u_i+hc/2 )\\ \hspace{-24pt} y_j=\phi(z_2,v_j)}  }.
\end{align*}
The property \eqref{associativitymod0} now follows by    examining the coefficients of the monomials \eqref{monomi} in \eqref{last4}. Indeed,  by the form of $\T^{  12}_{\bar{n}\bar{m}} $  and $\R^{{\text{\normalfont opp}}\ts 12}_{\bar{n}\bar{m}} $  there exists an integer $r_1\geqslant 1$ such that the coefficients of the given monomials in the product    $(z_1-z_2)^{r_1} \T^{  12}_{\bar{n}\bar{m}} \R^{{\text{\normalfont opp}}\ts 12}_{nm}$   possess only finitely many 
negative powers of the variables $z_1$ and $z_2$. Next, by the first family of relations in Proposition \ref{propositionz2} there exists an integer $r_2\geqslant 1$ such that the coefficients of the given monomials in the product    
$$(z_1-z_2)^{r_2}\ts \K_{[n]}^{- 13}(x)^{-1}\big|_{x_i=\phi(z_1,u_i+hc/2)}\ts ((RP)_{\bar{n}\bar{m}}^{12})^{-1} \ts
	\K_{[m]}^{-23}(y)^{-1}\big|_{y_j=\phi(z_2,v_j+hc/2) }$$
possess only finitely many 
negative powers of   $z_1$ and $z_2$. As the remaining terms in \eqref{last4} contain only nonnegative powers of $z_1$ and $z_2$, we conclude that the property \eqref{associativitymod0} holds for the monomials in \eqref{monomi} for any integer $r\geqslant r_1 + r_2$.

Finally, the weak associativity property \eqref{associativitymod} is   proved by similar calculations as the corresponding property from Definition \ref{mod}. However, in this case, its proof also   requires the use of explicit connection between the additive and multiplicative    $R$-matrices, \eqref{rat2} and \eqref{trig1}, as given in \eqref{conrmat} and \eqref{conrmats}. Note that the former $R$-matrix governs the $h$-adic quantum vertex algebra structure and the latter the $\phi$-coordinated module structure.\footnote{At this point, it is worth to recall that, in contrast,  the $h$-adic quantum vertex algebra structure and the corresponding module structure from Theorem \ref{glavni1} are both given in terms of the same  $R$-matrix  \eqref{rat2} with additive spectral parameter .} In fact, exactly because of the form of the substitution which appears on the left hand side of the identities \eqref{conrmat} and \eqref{conrmats},   the particular choice of  associate \eqref{assc5} enables us to establish the $\phi$-coordinated $\Vcr$-module structure on $\Mcrphi$.
\end{prf}

In the following theorem, we construct an example of $\phi$-coordinated $\Vcr$-module using the   module for the Yangian of RTT-type  established in Subsection \ref{sbsc42}. 
We omit its proof as it goes in parallel with the proofs of Theorems \ref{glavni1} and \ref{thm411}.
To emphasize that $\Ypbar
_c$ is now regarded as a $\phi$-coordinated $\Vcr$-module, we denote it by $\Mcrphir$. 

\begin{thm}\label{mcrrr}
For any $c\in\CC$ there exists a unique structure of $\phi$-coordinated module for the $h$-adic quantum vertex algebra $\Vcr$ on $\Mcrphir= \Ypbar_c$ such that for all integers $n\geqslant 1 $ and the variables $u=(u_1,\ldots ,u_n)$, $x=(x_1,\ldots ,x_n)$ we have
\beq\label{modmap3}
Y_{\Mcrphir}(T_{[n]}^+(u)\vac,z)
=
\ttt_{[n]}^+(x)\big|_{x_i=\phi(z,u_i)}\ts
\ttt_{[n]}^-(x)^{-1}\big|_{x_i=\phi(z,u_i+hc/2)}.
\eeq
\end{thm}

\subsection{On the algebra    \texorpdfstring{$\Dcrbar$}{D(R)c}}\label{subsec432}
In this subsection, we define a certain algebra via    relations which closely resemble those from the $RTT$-presentation of the quantum affine algebra in type $A$; see, e.g., \cite{FR,RS}.
For any $c\in\CC$ let  $\Dcrbar$ be the topologically free algebra over   $\CC[[h]]$ generated by  the elements $\bar{t}_{ij}^{(r)}$ with $i,j=1,\ldots ,N$ and $r\in\mathbb{Z}$. 
Its defining relations are  \eqref{rrt0}, \eqref{rrt} and \eqref{rrt2}, 
where 
$\wvr{T}^{+}(x)$ is defined as in \eqref{ghj76} and $\wvr{T}^{-}(x)$
is given by
\beq\label{6z7}
\wvr{T}^{-}(x)=\sum_{i,j=1}^N e_{ij}\ot \wvr{t}_{ij}^{-  }(x) \quad\text{for}\quad  
\wvr{t}_{ij}^-(x)= \sum_{r\geqslant 0} \wvr{t}_{ij}^{(r)} x^{-r}.
\eeq

The $\Dcrbar$-module $W$ is said to be {\em restricted} if it is topologically free as a $\CC[[h]]$-module   and the action of $\wvr{T}^-(x)$ on $W$  is invertible and such that  
$\wvr{T}^-(x)^{\pm 1}$ belongs to $\ndo\CC^N \ot \om(W,W[x^{-1}]_h)$. 
 Clearly, Proposition \ref{proorp} implies
\begin{kor}\label{korr2}
There exists a unique structure of restricted $\Dcrbar$-module
on $\Ypbar_c$   such that the action of $\wvr{T}^{-}(x)$ is given by \eqref{tubar} and   $\wvr{T}^+(x)$ acts by multiplication.
\end{kor}

Using   \eqref{rrt} one easily checks that on any restricted $\Dcrbar$-module $W$ we have 
$$
\wvr{T}_{[n]}^-(x_1,\ldots ,x_n)^{\pm 1}\in\ndo\CC^N \ot \om(W,W[x_1^{-1},\ldots ,x_{n}^{-1}]_h)\quad\text{for all }n\geqslant 1.
$$
Recall that  by  $\phi=\phi(z_2,z_0)$ we  denote  the associate given by \eqref{assc5}.  We have

\begin{kor}\label{sporedni2}
Let $W$ be a restricted $\Dcrbar$-module.  There exists a unique structure of $\phi$-coordinated   $\Vcr$-module on $W$ such that
\beq\label{action96}
Y_W(T_{[n]}^+(u)\vac,z)= \wvr{T}_{[n]}^+(x)\big|_{x_i=\phi(z,u_i )}\ts \wvr{T}_{[n]}^-(x)^{-1}\big|_{x_i=\phi(z,u_i+hc/2)}.
\eeq
Furthermore, if $ W_1\subset W$ is   a restricted $\Dcrbar$-submodule, then it is also a $\phi$-coordinated  $\Vcr$-submodule with respect to the suitable (co)restriction of $Y_W(\cdot,z)$.

\end{kor}

\begin{prf}
The first assertion 
can be proved by the same arguments as Theorem \ref{mcrrr}. More specifically, using the  defining relations  for the algebra $\Dcrbar$ and the fact that $W$ is restricted
one can   directly verify  the requirements from Definition \ref{phimod}. As for the second assertion, it is sufficient to observe that, due to  \eqref{action96},      for any $\CC[[h]]$-submodule $W_1\subset W$  the inclusions $\wvr{T}^{\pm}(x)W_1\subset \ndo\CC^N \ot W_1 [[x^{\pm 1}]]$   imply that
$Y_W(v,z)w$ belongs to $W_1[[z^{\pm 1}]]$ for all $v\in\Vcr$ and $w\in W_1$.
\end{prf}

By Corollary \ref{korr2}, $\Ypbar_c=\Mcrphir$ is a restricted $\Dcrbar$-module, so Corollary \ref{sporedni2} implies that it possesses a structure of $\phi$-coordinated $\Vcr$-module. However, in this  particular case, the corresponding $\phi$-coordinated module structure  coincides with the one which is already established by Theorem \ref{mcrrr}. 

\begin{rem}
The braided Yangian counterpart of   Corollary  \ref{sporedni2} can  be obtained as follows.
Using the family of  relations     \eqref{rttbraidedtrig}, \eqref{ll2bar} and \eqref{rll3bar}  one can again
introduce  an algebra over $\CC[[h]]$. Then, by arguing as in the proof  of   
 Theorem \ref{thm411}, one can show that all (suitably defined) restricted modules for the aforementioned algebra are naturally equipped with the structure of  $\phi$-coordinated  $\Vcr$-module.
It might be an interesting problem to investigate whether such algebras, along with braided Yangians,   can be realized as subalgebras of   $\Dcrbar$, or of the corresponding Yangians of the RTT-type, as is the case with certain classes of reflection algebras; cf. \cite{MR,S}. 
\end{rem}

\begin{rem}
If $R$ is a skew-invertible involutive symmetry, one can analogously use  relations \eqref{rtt}, \eqref{rtt2} and \eqref{rtt3}, given in terms of  the $R$-matrix \eqref{skiis}, to define a certain   algebra which resembles the double Yangian in type $A$; cf. \cite{I,Kho}. The   analogue of  Corollary \ref{sporedni2} for such algebra states that its  (suitably defined) restricted  modules are equipped with the structure of module for $\Vcr$ (in accordance with Remark \ref{skiisrem}). Thus, in particular, the $\phi$-coordinated modules are no longer required in such setting.
\end{rem}

\section{Quantum determinants}\label{sec732}

In this section, motivated by 
the explicit expressions for 
 the  determinants for the Yangians of   RTT-type and  the braided Yangians from \cite{GS}, we   obtain certain families of central elements  for  the $h$-adic quantum vertex algebra $\Vcr$ and  invariants for  its ($\phi$-coordinated)  modules from Sections \ref{subsec31} and \ref{subsec32}.

\subsection{Preliminaries on skew-symmetrizers}\label{sec732a}
Here we recall some properties of the skew-symmetrizer associated with the skew-invertible Hecke symmetry $R$ which we need later on. Our exposition closely follows \cite[Sect. 6]{GS}, which we adapt to our setting by taking $q=e^h$. To  simplify the notation, for any positive integer $k$ we shall write
$$
[k]_h = \frac{e^{kh} -e^{-kh}}{e^{h}-e^{-h}}.
$$
Note that $[k]_h$ belongs to $k+h\CC[[h]]$, so that we have
$[k]_h^{-1}\in \frac{1}{k}+h\CC[[h]]$.
 Define
\beq\label{rthat6}
\Rt(x)=R+\frac{(e^h - e^{-h})x}{1-x}I\in\ndo\CC^N\ot\ndo\CC^N[[x,h]].
\eeq
Clearly, $\Rt(x)$ is obtained from the $R$-matrix \eqref{trig1} by removing its normalization term and multiplying it by the permutation operator $P$ from the left. For any positive integer $m$ the {\em skew-symmetrizer} $\Pc^{(m)}$ associated to $R$ can be  defined as an element of $(\ndo\CC^N)^{\ot m}[[h]]$ such that
\beq\label{gss5}
\Pc^{(1)}=I\fand
\Pc^{(k+1)}= (-1)^k [k+1]_h^{-1} \Rt_{12}(e^{2h})\ts \Rt_{23}(e^{4h})
\ldots \Rt_{k\ts k+1}(e^{2kh})\ts \Pc^{(k)}_{1\ldots k}.
\eeq
On the right hand side, we use the usual notation convention, where the subscripts of elements indicate the tensor copies on which they are applied. For example, $\Pc^{(k)}_{1\ldots k}$ means that $\Pc^{(k)}$ acts on the first $k$ tensor factors of  $(\ndo\CC^N)^{\ot (k+1)}[[h]]$.

From now on, we assume that the rank of the skew-invertible Hecke symmetry $R$ is $(M|0)$; cf. \cite[Def. 2]{GPS}. 
Recall that the value of $M$ uniquely determines the normalization series from Subsection \ref{subsec1002}. 
By \cite[Lemma 19]{GS} we have  
\beq\label{ppom1}
\Pc^{(M)}_{2\ldots M+1}\ts \Rt_{1\to M}(x)
=(-1)^{M-1}\ts e^h\ts  [M]_h \ts \frac{1-xe^{-2Mh}}{ 1-xe^{-2(M-1)h}}\ts \Pc^{(M)}_{2\ldots M+1}\ts\Pc^{(M)}_{1\ldots M},
\eeq
where $\Rt_{1\to M}(x)$ stands for the product
$$
\Rt_{1\to M}(x)=\Rt_{12}(x)\ts \Rt_{23}(xe^{-2h})\ldots \Rt_{M\ts M+1}(xe^{-2(M-1)h}).
$$
Finally, due to \cite[Lemma 2.1]{HIOPT}, the skew-symmetrizer possesses the following properties:
\begin{gather*}
\Pc^{(M)}_{1\ldots M}\ts R_{M\ts M+1}\ldots R_{23}\ts R_{12}
=(-1)^{M-1}\ts e^h\ts [M]_h\ts 
\Pc^{(M)}_{1\ldots M}\ts
\Pc^{(M)}_{2\ldots M+1},\\
\Pc^{(M)}_{2\ldots M+1}\ts \Pc^{(M)}_{1\ldots M}\ts\Pc^{(M)}_{2\ldots M+1}
=[M]_h^{-2}\ts \Pc^{(M)}_{2\ldots M+1}.
\end{gather*}

As $R$ is a skew-invertible Hecke symmetry of rank $(M|0)$, the image of the  skew-symmetrizer
$\Pc^{(M)}$ is spanned by a single nonzero vector. Let us denote one such vector by $v_R$. It can be proved that there exists a linear map $\mathcal{N}_R$ such that for all $w\in \CC^N$ we have
$
R_{12}\ts R_{23}\ldots R_{M\ts M+1} (v_R \ot w) = \mathcal{N}_R (w)\ot v_R
$; see \cite{Gur,GS} for more details and the explicit expression for 	 $\mathcal{N}_R$.
 
\subsection{Center of $\Vcr$}\label{sec732b}

Recall that by $R$ we denote  a skew-invertible Hecke-symmetry   of rank $(M|0)$ and by $\Psi$ its skew-inverse  \eqref{skew}. Let $C=\tr_2 \Psi$, where the subscript $2$ indicates that the trace is taken over the second tensor factor. Thus, $C=(c_{ij})$ is the square matrix of order $N$.
The {\em $R$-trace} $\tr^R A$ of an arbitrary $N\times N$ matrix $A=(a_{ij})$ is defined by
$$
\tr^R A= \tr \left(A\cdot C\right)=\sum_{i,j=1}^N a_{ij}c_{ji}.
$$ 
It   is generalized to   multiple tensor factors by
setting $\tr_{1,\ldots,k}^R =\tr_{1}^R\ldots \tr_{k}^R$, where $\tr_i^R$ denotes the $R$-trace taken over the $i$-th tensor factor.
Motivated by the form of the  determinant for the Yangian of   RTT-type \cite[Rem. 22]{GS}, we consider the   power series
$$
\qdet T^+(u)=\tr_{1,\ldots,M}^R\ts \Pc^{(M)} \ts T_1^+(u)\ts T_2^+(u+ah)\ldots T_M^+(u+a(M-1)h ) \in\ndo \Vcr[[u]].
$$
Note that the given expression is well-defined because $\Vcr$ is $h$-adically complete.

Before we state our next result, we recall that the {\em center} $\z(V)$ 
of an  $h$-adic quantum vertex algebra $V$  (see \cite{DGK,JKMY} for more information) is defined by
$$
\z(V)=\left\{w\in V\, :\, Y(v,z)w\in V[[z]]\text{ for all }v\in V\right\}.
$$

\begin{thm}\label{centralni}
If $\mathcal{N}_R$ is a scalar matrix, then
all coefficients of the series $\qdet T^+(u)\vac$ belong to the center $\z(\Vcr)$ for any $c\in\CC$.
\end{thm}

\begin{prf}
{\em (1)} Suppose that the following identity for operators on $\Vcr$ holds:
\beq\label{gss2}
T^-(u)\ts \qdet T^+(v) =\qdet T^+(v)\ts T^-(u).
\eeq
It implies that for any integer $n\geqslant 1$ and the variables $u=(u_1,\ldots ,u_n)$ and $z$ we have
$$T_{[n]}^-(z+u+hc/2)^{-1}\ts \qdet T^+(v)
=\qdet T^+(v)\ts T_{[n]}^-(z+u+hc/2)^{-1}.$$
By applying this on the vacuum vector $\vac$ and using $T^-(u)^{-1}\vac=I\ot \vac$ we get
$$T_{[n]}^-(z+u+hc/2)^{-1}\ts \qdet T^+(v)\vac
=\qdet T^+(v)\vac.$$
Next, using the explicit expression \eqref{vmap} for the vertex operator map we find
\beq\label{gss1}
Y(T_{[n]}^+(u),z)\ts \qdet T^+(v)\vac =T_{[n]}^+(z+u)\ts\qdet T^+(v)\vac.
\eeq
The right hand side contains only nonnegative powers of the variable $z$. Moreover, 
as the coefficients of all matrix entries of $T_{[n]}^+(u_1,\ldots ,u_n)$ with $n\geqslant 1$ along with $\vac$ span an $h$-adically dense $\CC[[h]]$-submodule of $\Vcr$, we conclude from \eqref{gss1} that
$Y(w,z)\qdet T^+(v)\vac$ belongs to $\Vcr[[z]]$ for all $w\in\Vcr$, so that all coefficients of $\qdet T^+(v)\vac$ are elements of the center $\z(\Vcr)$. Therefore, in order to prove the theorem, it is sufficient to check that the equality \eqref{gss2} holds. This can be done by arguing as in the proof of \cite[Prop. 21]{GS}; see also \cite[Rem. 22]{GS}. However, we present   the underlying calculations below in order to take care of   differences which occur in our setting, such as the dependency on the properties of the normalization series from Subsection \ref{subsec1002} (when $c$ is nonzero).

{\em \noindent(2)} Conjugating  relation \eqref{rtt3} by the permutation operator $P$ and using \eqref{U} we get
\begin{align*}
&f(-u+v-hc/2)\ts
\Rt_{12}(e^{-\frac{2}{a}(-u+v-hc/2)})\ts T_1^+(v)\ts T_2^- (u)\\
&\qquad=
f(-u+v+hc/2)\ts
T_1^- (u)\ts T_2^+(v)\ts \Rt_{12}(e^{-\frac{2}{a}(-u+v+hc/2)}).
\end{align*}
From this, one easily derives the more general identity
\begin{align}
&F^-\ts \Rt_{1\to M}(e^{-\frac{2}{a}(-u+v-hc/2)})\ts T_1^+(v)\ts T_2^+(v+ah)\ldots T_M^+(v+a(M-1)h)\ts T_{M+1}^-(u)\non\\
=& \, T_{1}^-(u)\ts T_2^+(v)\ts T_3^+(v+ah)\ldots T_{M+1}^+(v+a(M-1)h)\ts  \Rt_{1\to M}(e^{-\frac{2}{a}(-u+v+hc/2}))\ts F^+,\label{gss3}
\end{align}
where, due to the property  \eqref{normal5} of the series $f(u)$, the expressions $F^{\pm}$ are given by
\beq\label{normal8}
F^\pm =\prod_{i=1}^M f (-u+v+a(i-1)h\pm hc/2)=\frac{1-e^{-\frac{2}{a}(-u+v\pm hc/2)} e^{-2(M-1)h}}{1-e^{-\frac{2}{a}(-u+v\pm hc/2)} e^{-2M h}}.
\eeq
By combining \eqref{ppom1} and \eqref{normal8} we find
\begin{align}
&\Pc^{(M)}_{2\ldots M+1}\ts F^\pm\ts \Rt_{1\to M}(e^{-\frac{2}{a}(-u+v\pm hc/2)}) =(-1)^{M-1}\ts e^h\ts  [M]_h \ts   \Pc^{(M)}_{2\ldots M+1}\ts\Pc^{(M)}_{1\ldots M}.\label{gss7}
\end{align}
Therefore, multiplying \eqref{gss3} by $\Pc^{(M)}_{2\ldots M+1}$ from the left gives us
\begin{align}
&(-1)^{M-1}\ts e^h\ts  [M]_h\ts \Pc^{(M)}_{2\ldots M+1}\ts\Pc^{(M)}_{1\ldots M}\ts T_1^+(v) \ldots T_M^+(v+a(M-1)h)\ts T_{M+1}^-(u)\non\\
=&\, T_{1}^-(u)\ts \Pc^{(M)}_{2\ldots M+1}\ts T_2^+(v) \ldots T_{M+1}^+(v+a(M-1)h)\ts  \Rt_{1\to M}(e^{-\frac{2}{a}(-u+v+hc/2}))\ts F^+.\label{gss6}
\end{align}
Next, by using the identities \eqref{rtt} and \eqref{gss5} we obtain
\begin{align*}
&\Pc^{(M)}_{1\ldots M}\ts T_1^+(v)\ts T_2^+(v+ah)\ldots T_{M}^+(v+a(M-1)h)\\
&\qquad = T_{1}^+(v+a(M-1)h)\ldots T_{M-1}^+(v+ah)\ts T_M^+(v)\ts \Pc^{(M)}_{1\ldots M}.
\end{align*}
Applying this to the tensor factors $2,\ldots, M+1$ on the right hand side of  \eqref{gss6}, then using \eqref{gss7} and, finally,  canceling    common factors on both sides,   \eqref{gss6} takes the form  
\begin{align*}
&  \Pc^{(M)}_{2\ldots M+1}\ts\Pc^{(M)}_{1\ldots M}\ts T_1^+(v)\ts T_2^+(v+ah)\ldots T_M^+(v+a(M-1)h)\ts T_{M+1}^-(u)\non\\
&\qquad = T_{1}^-(u)\ts  T_2^+(v+a(M-1)h)\ts  \ldots T_M^+(v+ah)\ts T_{M+1}^+(v)\ts   \Pc^{(M)}_{2\ldots M+1}\ts\Pc^{(M)}_{1\ldots M}. 
\end{align*}
As in the corresponding part of the proof of \cite[Prop. 21]{GS}, this turns to the equality
$$
 \Pc^{(M)}_{2\ldots M+1}\ts\Pc^{(M)}_{1\ldots M}\ts \qdet  T^+(v)\ts T_{M+1}^-(u)\ = T_{1}^-(u)\ts  \ts \qdet  T^+(v)\ts   \Pc^{(M)}_{2\ldots M+1}\ts\Pc^{(M)}_{1\ldots M	}, 
$$
which can be transformed as in   \cite[Rem. 22]{GS} into 
$$
 \mathcal{N}_R\ts T^-(u)\ts\qdet  T^+(v)  = \qdet  T^+(v)\ts   T^-(u)\ts   \mathcal{N}_R. 
$$
Hence, if the matrix $\mathcal{N}_R$ is scalar, this implies \eqref{gss2}, as required.
\end{prf}

\subsection{Invariants of (\texorpdfstring{$\phi$}{phi}-coordinated) \texorpdfstring{$\Vcr$}{Vc(R)}-modules}\label{sec732c}
In this subsection, we consider the image of the constant term of the quantum determinant under the ($\phi$-coordinated) $\Vcr$-module maps \eqref{modmap},  \eqref{modmap2} and \eqref{modmap3}. 
We shall write $(RP)_{[1]}=I\in \ndo\CC^N$ and
$$
(RP)_{[M]} =\prod_{1\leqslant i< M }^{\longrightarrow}
  (RP)_{i\ts i+1}\ldots (RP)_{i\ts M}\in(\ndo\CC^N)^{\ot M}[[h]]\qquad\text{for}\qquad M>1,
$$
where, as before,    $R$ is a skew-invertible Hecke-symmetry  of rank $(M|0)$ and $I $   the identity matrix.
Notice that for   $u=(u_1,\ldots,u_M)$ and $x=(x_1,\ldots,x_M)$ we have
\beq\label{obs}
L_{[M]}^- (u )^{-1}1 = \K_{[M]}^- (x )^{-1}1
=(RP)_{[M]}^{-1}\ot 1
\Fand
\wvr{T}^-_{[M]} (x )^{-1}1=I^{\ot M} \ot 1
,
\eeq
where $I\in \ndo\CC^N$ denotes the identity matrix of order $N$. We now employ the aforementioned ($\phi$-coordinated) $\Vcr$-module maps to introduce the following power series
\begin{align*}
&\qdet L^+(z)= Y_{\Mcr}(\qdet T^+(0)\vac,z)1\in \Mcr[[z]],\\
&\qdet \K^+(z)= Y_{\Mcrphi}(\qdet T^+(0)\vac,z)1\in \Mcrphi[[z]],\\
&\qdet \ttt^+(z)= Y_{\Mcrphir}(\qdet T^+(0)\vac,z)1\in \Mcrphir[[z]]. 
\end{align*}
As   suggested by our   notation, we shall refer to them as quantum determinants. 

By using  the expressions for the  ($\phi$-coordinated) $\Vcr$-module maps \eqref{modmap},  \eqref{modmap2} and \eqref{modmap3} and the identities in \eqref{obs} one easily derives explicit formulae
\begin{align}
&\qdet L^+(z)=\tr_{1,\ldots,M}^R\ts \Pc^{(M)} \ts L_{[M]}^+(z,z+ah,\ldots ,z+a(M-1)h )\ts (RP)_{[M]}^{-1},\label{allignnn}\\
&\qdet \K^+(z)=\tr_{1,\ldots,M}^R\ts \Pc^{(M)} \ts \K_{[M]}^+(z,ze^{-2h},\ldots ,ze^{-2(M-1)h} )\ts (RP)_{[M]}^{-1},\label{allignn}\\
&\qdet \ttt^+(z)=\tr_{1,\ldots,M}^R\ts \Pc^{(M)} \ts \ttt_{[M]}^+(z,ze^{-2h},\ldots ,ze^{-2(M-1)h} ) . \label{allign}
\end{align}
Therefore, the  determinants \eqref{allignnn}, \eqref{allignn} and \eqref{allign} can be naturally regarded as operators on the $\CC[[h]]$-modules  $\Mcr$, $\Mcrphi$  and  $\Mcrphir$ with respect to the multiplication in the corresponding algebra  $\Ybp$, $\Ybpbar$ and $\Ypbar$.

\begin{lem}\label{yg}
Let $c\in\CC$.
The next identities hold for operators on $\Mcr$ and $\Mcrphi$:
$$
\qdet L^+(z_1)\ts L^-(z_2)=L^-(z_2)\ts \qdet L^+(z_1)\fand
\qdet \K^+(z_1)\ts \K^-(z_2)=\K^-(z_2)\ts \qdet \K^+(z_1).
$$
If $\mathcal{N}_R$ is a scalar matrix, we have the following identity for operators on $\Mcrphir$:
$$
\qdet \ttt^+(z_1)\ts \ttt^-(z_2)=\ttt^-(z_2)\ts \qdet \ttt^+(z_1).
$$
\end{lem}

\begin{prf}
The   lemma is verified by the arguments which closely follow the proofs of \cite[Prop. 21, Rem. 22]{GS} and   Theorem \ref{centralni} and rely on the properties of the normalizing series  and the skew-symmetrizer from Subsections  \ref{subsec1002}   and    \ref{sec732a}.
In order to employ these arguments, it is useful to observe that, despite seemingly different expressions, 
  $\qdet \K^+(z)$   is of the same form as the series    $e_M (z)$ which belongs to the center of the braided Yangian \cite[Prop. 21]{GS}. Indeed,   this is easily  proved by moving all   permutation operators $P$ which appear in the above expression  for 
 $\qdet \K^+(z)$   to the right and then using the braid relation \eqref{hecke1}, which is satisfied by $R$.  The same remark applies to the other two determinants.
\end{prf}

We  are now  ready to present the vertex algebraic interpretation of quantum determinants.
Suppose $W$ is a ($\phi$-coordinated)  module for the $h$-adic quantum vertex algebra $V$ with respect to the map $Y_W(\cdot ,z )$. Define the {\em submodule of invariants} $\z(W)$ of $W$ by
$$
\z(W)=\left\{w\in W\,:\, Y_W(v,z)w\in W[[z]]\text{ for all }v\in V\right\}.
$$
The following theorem shows that the quantum determinants produce invariants of ($\phi$-coordinated) $\Vcr$-modules.
Due to Lemma \ref{yg}, it can be verified  by   arguing as in the first part of the proof of Theorem \ref{centralni}.

\begin{thm}
Let $c\in\CC$. All coefficients of   $\qdet L^+(z)1$ and $\qdet\K^+(z)1$ belong to the corresponding submodule of invariants
$\z(\Mcr)$ and $\z(\Mcrphi)$, respectively. Moreover, if $\mathcal{N}_R$ is a scalar matrix, then all coefficients of $\qdet \ttt^+(z)1$ belong to $\z(\Mcrphir)$.
\end{thm}

Due to \eqref{allign} the series  $\qdet \ttt^+ (x)$ can be also regarded as a  power series with coefficients in $\Dcrbar$, i.e. as an element of $\Dcrbar[[x]]$. Thus,  we have the following corollary:  

\begin{kor}
 If $\mathcal{N}_R$ is a scalar matrix, then all coefficients of $\qdet \ttt^+(x)$ belong to the center of the algebra $\Dcrbar$ for any $c\in\CC$.
\end{kor} 

\begin{prf}
Clearly, it  suffices to show that $\qdet \ttt^+(x)$ commutes with the generator matrices $\ttt^-(y)$ and $\ttt^+(y)$ of $\Dcrbar$.
 Regarding the former matrix, the commutation relations    between $\ttt^+(x)$ and $\ttt^-(y)$, when regarded as operators on  $\Mcrphir$ and when regarded as  series  with coefficients in $\Dcrbar$, coincide, so this follows  from  Lemma \ref{yg}.
As for the latter matrix, the family  of defining relation \eqref{rrt0}    can be    written in terms of \eqref{rthat6} as
\beq\label{thsdf}
\Rt_{12}(x/y)\ts \ttt_1^+(x)\ts   \ttt_2^+(y)=\ttt_1^+(y) \ts \ttt_2^+(x) \ts \Rt_{12}(x/y). 
\eeq
Indeed, the above equality is found  
by canceling the normalization series $\wvr{f}(x/y)$ in \eqref{rrt0} and then multiplying the resulting equality by the permutation operator   from the left. As the form of \eqref{thsdf} coincides with the form of the defining relation for the Yangian of   RTT-type, one   shows that $\qdet \ttt^+(x)$ and 
$\ttt^+(y)$ commute     as in \cite[Rem. 22]{GS}.
\end{prf}

\section*{Acknowledgement}
This work has been supported in part by Croatian Science Foundation under the project UIP-2019-04-8488.

\end{document}